\documentstyle[12pt]{article}
\begin{document}

\def\dsp{\displaystyle}
\def\Rr{{bf R}} 
\def\Zz{bf Z}
\def\Nn{bf N}
\def\get{\hbox{{\goth g}$^*$}}
\def\g{\gamma}
\def\om{\omega}
\def\r{\rho}
\def\a{\alpha}
\def\s{\sigma}
\def\vfi{\varphi}
\def\l{\lambda}
\def\implique{\Rightarrow}
\def\o{{\circ}}
\def\Diff{\hbox{\rm Diff}}
\def\S1{\hbox{\rm S$^1$}}
\def\Hom{\hbox{\rm Hom}}
\def\Vect{\hbox{\rm Vect}}
\def\const{\hbox{\rm const}}
\def\ad{\hbox{\hbox{\rm ad}}}
\def\semid{\hbox{\bb o}}
\def\blanc{\hbox{\ \ }}

\def\pds#1,#2{\langle #1\mid #2\rangle} 
\def\f#1,#2,#3{#1\colon#2\to#3} 

\def\hfl#1{{\buildrel{#1}\over{\hbox to
12mm{\rightarrowfill}}}}

\title{Projectively invariant symbol map and
cohomology of vector fields Lie algebras intervening in quantization}

\author{P.B.A. Lecomte
\thanks{Institute de Math\'ematiques,
 Universit\'e de Li\`ege,
Avenue des Tilleuls, 15,
B-400 Li\`ege,
 BELGIQUE
}
\and
 V.Yu. Ovsienko
\thanks{C.N.R.S., Centre de Physique Th\'eorique,
 Luminy -- Case 907,
F--13288 Marseille, Cedex 9,
FRANCE
}
}

\date{}

\maketitle

{\abstract{
We define the unique (up to normalization) 
symbol map from the space of
linear differential operators on ${\bf R}^n$
to the space of polynomial on fibers functions on $T^*{\bf R}^n$,
equivariant with respect to the Lie algebra
of projective transformations 
$sl_{n+1}\subset\Vect({\bf R}^n)$.
We apply the constructed $sl_{n+1}$-invariant
symbol to studying of the
natural one-parameter
family of $\Vect(M)$-modules on the space of linear differential operators
on an arbitrary manifold $M$.
Each of the $\Vect(M)$-action from this family
can be interpreted as a deformation of
the standard $\Vect(M)$-module ${\cal S}(M)$
 of symmetric contravariant tensor fields on $M$.
We define (and calculate
in the case: $M={\bf R}^n$)
the corresponding cohomology of $\Vect(M)$
related with this deformation.
This cohomology realize the obstruction for existence of equivariant 
symbol and quantization maps.
The projective Lie algebra  $sl_{n+1}$
naturally appears as the algebra of symmetries
on which the involved $\Vect(M)$-cohomology is trivial.
}}

\vfill\eject

\oddsidemargin  -10pt    
\hsize 6.4truein  
\voffset=-1truecm

\setlength{\textheight}{21cm}

{\Large{\bf Introduction}}

\vskip 0,5cm

The general problem of quantization
is
to associate 
linear operators on a Hilbert space
with Hamiltonian functions
on a symplectic manifold.
In the simplest situation 
when the symplectic manifold is a cotangent bundle $T^*M$
over a manifold $M$,
the usual quantization procedure consists of establishing 
a correspondence between the space ${\rm Pol}(T^*M)$ of functions on $T^*M$
polynomial on fiber and the space ${\cal D}(M)$ 
of linear differential operators on $M$:
$$
\sigma: {\cal D}(M)\to {\rm Pol}(T^*M)
\;\;\;\;\;
\hbox{symbol map}.
$$
and the inverse
$$
\sigma^{-1}: {\rm Pol}(T^*M)\to{\cal D}(M)
\;\;\;\;\;
\hbox{quantization map}
$$
(see e.g. \cite{ber}).

\vskip 0,3cm

Let $\Diff(M)$ be the group of diffeomorphisms of $M$
and $\Vect(M)$ be the Lie algebra of vector fields on $M$.
Space ${\cal D}(M)$ has
a natural 1-parameter family of 
$\Diff(M)$- and $\Vect(M)$-modules.
To define it, one considers arguments of differential operators as tensor-densities
of an arbitrary degree $\lambda$ (see \cite{duv},\cite{lec},\cite{gar}).

One of the difficulties of quantization 
can be formulated as follows.
{\it There is no quantization map (and symbol map) equivariant with respect to
the action of} $\Diff(M)$.
Two spaces: ${\rm Pol}(T^*M)$ and ${\cal D}(M)$
are not isomorphic as modules over
$\Diff(M)$ or $\Vect(M)$.
To avoid this difficulty,
some of methods of quantization are based on
a choice of Darboux coordinates in $T^*M$
(e.i. the well-known Weyl quantization
and the Moyal-Weyl deformation quantization).
Other ones
(as the modern approach of B. Fedosov \cite{fed})
fix a linear (symplectic) connection on $T^*M$.
Another way is to use
the representation theory
(as in the case of geometric quantization).

\vskip 0,3cm

This paper contains two main parts.

\vskip 0,3cm

(a) 
We consider a natural embedding
$sl_{n+1}\subset\Vect({\bf R}^n)$
($sl_{n+1}$ acts on ${\bf R}^n$
by infinitesimal projective transformations).

The first main result of this paper is the existence of the unique
(up to normalization)
$sl_{n+1}$-equivariant symbol map
on ${\bf R}^n$.
This implies that ${\cal D}({\bf R}^n)$
is isomorphic to ${\rm Pol}(T^*{\bf R}^n)$
as a module over $sl_{n+1}$.

The same result holds for manifolds
endowed with a projective structure
(an atlas with linear-fractional coordinate transformations).
The main examples are:
$M={\bf P}^n,S^n,{\bf T}^n$.

\vskip 0,3cm

{\bf Remarks}. 1. In the one-dimensional case ($n=1$), 
$sl_2$-equivariant symbol map and quantization map
was obtained (in a more general situation of pseudodifferential operators)
in a recent work by P.B. Cohen, Yu. I. Manin and D. Zagier \cite{cmz}.
We check that in the case $n=1$, our formul{\ae} coincide with those of \cite{cmz}.
We were not aware of this paper while the calculation
of the $sl_{n+1}$-equivariant symbol has been done.

2. In the (algebraic) case of global differential operators
on ${\bf CP}^n$, existence and uniqueness of the 
$sl_{n+1}$-equivariant symbol is a corollary of 
Borho--Brylinski's results \cite{bor}.
Space ${\cal D}({\bf CP}^n)$ as a module over $sl_{n+1}$
has a decomposition in a sum of irreducible submodules
with {\it multiplicities one}. This implies the uniqueness result.
Our explicit formul{\ae} are valid in the holomorphic case and
define an isomorphism between ${\cal D}({\bf CP}^n)$
and the space of polynomial on fibers functions on $T^*{\bf CP}^n$.

\vskip 0,3cm

(b)
In the second part of this paper,
we study the one-parameter family of $\Vect(M)$-modules 
on ${\cal D}(M)$
for an arbitrary smooth manifold $M$.

Space ${\rm Pol}(T^*M)$ of fiber-vise polynomials on $T^*M$
is naturally isomorphic (as a $\Diff(M)$-module) to
the space of symmetric contravariant tensor fields on $M$:
$\Gamma({\rm S}(TM))$.
Action of $\Diff(M)$ and $\Vect(M)$ on ${\cal D}(M)$
can be realized as a {\it deformation} of the standard action on ${\rm Pol}(T^*M)$.
This view point makes possible to apply the cohomology
technique: 
the action of $\Vect(M)$ on ${\cal D}(M)$
is distinguished from the action of $\Vect(M)$ on 
${\rm Pol}(T^*M)$ by certain
cohomology classes of $\Vect(M)$.
More precisely, this approach leads to the first group 
of $\Vect(M)$-cohomology with
coefficients in the space of linear
operators on symmetric contravariant tensor fields.
This cohomology realizes
obstructions
for existence of equivariant quantization (and symbol) map.

The cohomological approach to studying of 
$\Vect(M)$-module structure on the space of differential operators
was proposed in \cite{duv} in the case of second order operators.

\vskip 0,3cm

The second main result of this paper is
calculation of the first group of
$\Vect({\bf R}^n)$-cohomology 
vanishing on $sl_{n+1}$,
with coefficients in the space
of differentiable linear operators on symmetric
contravariant tensor fields.
This result includes classification of
bilinear $sl_{n+1}$-equivariant maps
$\Vect({\bf R}^n)\otimes\Gamma(S^k(TM))\to\Gamma(S^l(TM))$
vanishing on subalgebra $sl_{n+1}\subset\Vect({\bf R}^n)$.
Such operations are multi-dimensional analogues
of so-called Gordan transvectants \cite{gor}
(bilinear $sl_2$-equivariant operators on tensor-densities on $S^1$).

\vskip 0,3cm

 Lie algebra $sl_{n+1}$ appears as
the maximal subalgebra of $\Vect({\bf R}^n)$ on which all the
involved cohomology classes vanish. The restriction to
 $sl_{n+1}$ of the 
$\Vect({\bf R}^n)$-modules on ${\cal D}({\bf R}^n)$,
is trivial: all the $sl_{n+1}$-modules on ${\cal D}({\bf R}^n)$
 are isomorphic to
the module of symmetric contravariant tensor fields.

\vskip 0,3cm

The problem of isomorphism of the $\Vect(M)$-modules
on ${\cal D}(M)$
for different values of the parameter $\lambda$ 
was solved in a series of recent papers
\cite{duv},\cite{lec},\cite{gar}.
We will classify the factor-modules
${\cal D}^k_{\lambda}(M)/{\cal D}^l_{\lambda}(M)$
(of the module of $k$-th order operators 
over the module of $l$-th order operators).

\vskip 0,3cm

(c) We will apply the $sl_{n+1}$-equivariant symbol and
quantization maps to 
define $sl_{n+1}$-equivariant quantization
of $T^*S^n$. Namely, we will define a $sl_{n+1}$-equivariant
star-product on $T^*S^n$
and give a $sl_{n+1}$-equivariant version of
quantization of the geodesic flow on $T^*S^n$.

\vskip 0,3cm

We hope, that the appearance of projective symmetries
in the context of quantization
is natural.
We do not know any works on this subject except the special
one-dimensional case (cf. \cite{cmz},\cite{ovs}).

Relations of differential
operators with projective geometry have already been studied by classics
(see e.g. \cite{wil},\cite{car}). The most known example is the
Sturm-Liouville operator $d^2/dx^2+u(x)$ describing a projective structure
on $\bf R$ (or on $S^1$ if $u(x)$ is periodic).
This example is related to the well-known Virasoro algebra (see \cite{kir}).

\section{Projective Lie algebra $sl_{n+1}$}

We call the {\it projective Lie algebra},
the Lie algebra of vector fields on ${\bf R}^n$
generated 
by the following vector fields:
\begin{equation}
\frac{\partial}{\partial x^i},
\;\;\;x^i\frac{\partial}{\partial x^j},
\;\;\;x^iE,
\label{sl}
\end{equation}
where
$$
E=x^j\frac{\partial}{\partial x^j}
$$ 
(we omit the sign of a sum over repeated indices)
is the Euler field.
This Lie algebra is isomorphic to $sl_{n+1}$.
The vector fields (\ref{sl}) 
are infinitesimal generators of (locally defined on ${\bf R}^n$)
linear-fractional transformations.

The constant and linear vector fields (first two terms in (\ref{sl}))
generate a subalgebra of $sl_{n+1}$ called the {\it affine Lie algebra}.

\vskip 0,3cm

The projective Lie algebra
is well-defined
globally on ${\bf P}^n\supset{\bf R}^n$.
The space of vector fields generated by the vector fields (\ref{sl})
is invariant under linear-fractional coordinate changes.

\vskip 0,3cm

{\bf Remark 1.1}.
Lie subalgebra $sl_{n+1}\subset\Vect({\bf R}^n)$ is {\it maximal} in the following sense.
Given an arbitrary 
polynomial vector field $X\not\in sl_{n+1}$,
then, a Lie algebra generated
by $X$ together with the vector fields (\ref{sl}),
coincides with the whole Lie algebra of polynomial
vector fields on ${\bf R}^n$
(see \cite{ogi} for a similar result).

\vskip 0,3cm

Lie algebra $sl_{n+1}$
can be defined as the kernel of
certain cohomology classes of $\Vect({\bf R}^n)$.
This interesting relation between 
the cohomology of $\Vect({\bf R}^n)$
and the projective Lie algebra $sl_{n+1}$ is
the main tool of this paper.

\subsection{$sl_{n+1}$ as a cohomology kernel}

Consider the space ${\cal S}^1_2$ 
of symmetric $(1,2)$-tensor fields on ${\bf R}^n$:
$$
T=T_{ij}^k(x)\cdot dx^idx^j\otimes \partial_k,
$$
here and below $\partial_k=\frac{\partial}{\partial x^k}$.
There exist a non-trivial cocycle
$\gamma:\Vect({\bf R}^n)\to{\cal S}^1_2$:
$$
\gamma(X)
=
\partial_i\partial_j(X^k)\cdot
dx^idx^j\otimes\partial_k,
$$
where
$X=X^k\frac{\partial}{\partial x^k}$.

Cocycle $\gamma$ can be written in the following form.
Let $D$ be the differential,
acting on symmetric tensor fields:
$D(a)=\partial_i(a)\cdot dx^i$
(operator $D$ depends on the choice of coordinates).
Then,
$
\gamma(X)=D^2(X),
$

\vskip 0,3cm

{\bf Remark 1.2}. Note, that cocycle $\gamma$ 
is just the non-linear term 
in the standard (infinitesimal)
coordinate transformations of the Christoffel symbols
of a linear connection.

\vskip 0,3cm

The important property of the projective subalgebra 
$sl_{n+1}\subset\Vect({\bf R}^n)$
is that the restriction to $sl_{n+1}$ of
the cohomology class $[\gamma]$
 vanishes.

\proclaim Proposition 1.3. (i) In the multi-dimensional case ($n\geq2$),
there exists a unique (up to a constant) cocycle
on $\Vect({\bf R}^n)$ with values in ${\cal S}^1_2$
vanishing on $sl_{n+1}$:
\begin{equation}
\bar\gamma
(X)=
\gamma(X)-
\frac{2}{n+1}\partial_i\partial_k(X^k)\cdot
dx^idx^j\otimes\partial_j
\label{ha}
\end{equation}
$\bar\gamma$ is nontrivial and cohomological to $\gamma$.
\hfill\break
(ii) The vector fields (\ref{sl}) generate the space
of solutions of the equation $\bar\gamma(X)=0$.
\par

{\bf Proof}: The uniqueness of $\bar\gamma$ is a corollary of Theorem 5.2 below.

\vskip 0,3cm

The cocycle $\bar\gamma$ is a projection of
$\gamma$ to the space of zero-trace tensor fields (such that $T_{ij}^i\equiv0$).
To rewrite $\bar\gamma$ in an intrinsic way, introduce
a $(1,1)$-tensor $\Delta=dx^j\otimes\partial_j$. Then,
$
\bar\gamma
(X)=
\gamma(X)-
\frac{2}{n+1}D({\rm div}X)\Delta .
$

\vskip 0,3cm

The cocycle (\ref{ha}) was considered in \cite{hab}.

\subsection{Projectively equivariant cocycles}

A general property of 1-cocycles is that
{\it a 1-cocycle vanishing on a Lie subalgebra 
 is equivariant with respect to this
subalgebra.}
In particular, the cocycle $\bar\gamma$
on $\Vect({\bf R}^n)$
defines a $sl_{n+1}$-equivariant map
from $\Vect({\bf R}^n)$ to ${\cal S}^1_2$.

\vskip 0,3cm

Indeed, the 1-cocycle relation:
$
L_X(\bar\gamma(Y))-L_Y(\bar\gamma(X))
=
\bar\gamma([X,Y])
$
implies the property of $sl_{n+1}$-equivariance:
$
L_X(\bar\gamma(Y))
=
\bar\gamma([X,Y])$,
for every $X\in sl_{n+1}$.

\subsection{Gelfand-Fuchs cocycle}

In the one-dimensional case,
$\bar\gamma\equiv0$ and
and Proposition 1.2 does not hold.
The cocycle $\gamma$ in this case is trivial.

However, there exists a non-trivial
cocycle on $\Vect({\bf R})$
with values in the quadratic differentials
vanishing on the subalgebra $sl_2$:
\begin{equation}
X\frac{d}{dx}\mapsto X'''\;(dx)^2.
\label{gf}
\end{equation}
This is (a version of) so-called {\it Gelfand-Fuchs cocycle} (see e.g. \cite{fuc})
related to the Virasoro algebra.

We will define a multi-dimensional analogue of the
cocycle (\ref{gf}) in Section 5.4.

\section{Equivariant symbol}

Let ${\cal D}$ be the space of
scalar linear differential operators
on ${\bf R}^n$ :
\begin{equation}
A
=
a_k^{{i_1}\ldots{i_k}}
\frac{\partial}{\partial x^{i_1}}\cdots
\frac{\partial}{\partial x^{i_k}}
+
\cdots
+
a_1^i
\frac{\partial}{\partial x^i}
+
a_0
\label{A}
\end{equation}
with coefficients 
$a_j^{{i_1}\ldots{i_j}}=a_j^{{i_1}\ldots{i_j}}(x^1,\dots,x^n)\in C^{\infty}({\bf R}^n)$,
number $k$ is the order of $A$.

Denote ${\cal D}^k\subset{\cal D}$ the space of all $k$-th order
linear differential operators.

\subsection{$\Diff({\bf R}^n)$- and $\Vect({\bf R}^n)$-module structures}

Let us recall the definition of
the natural 1-parameter family of 
$\Diff({\bf R}^n)$- and $\Vect({\bf R}^n)$-module structures
on ${\cal D}$ (see \cite{duv},\cite{lec},\cite{gar}).

\vskip 0,3cm

Consider a 1-parameter family 
of $\Diff({\bf R}^n)$-actions on $C^{\infty}({\bf R}^n)$:

$$
g^*_{\lambda}(\phi)
:=\phi\circ g^{-1}\cdot\left|\frac{Dg^{-1}}{Dx}\right|^{\lambda},
$$
where $g\in\Diff({\bf R}^n)$, $\phi\in C^{\infty}({\bf R}^n)$,
$\left|Dg^{-1}/Dx\right|$ is the Jacobian of $g^{-1}$ and
$\lambda\in {\bf R}$ (or ${\bf C}$).

The 1-parameter family of $\Diff({\bf R}^n)$-actions on ${\cal D}$ is given by:
\begin{equation}
g^{\lambda}(A)=g^*_{\lambda}A(g^*_{\lambda})^{-1}.
\label{act}
\end{equation}

\vskip 0,3cm

{\bf Remark 2.1}. Operator $g^*_{\lambda}$ is
the natural action of $\Diff({\bf R}^n)$ on the space of
{\it tensor-densities} of degree
 $\lambda$ on ${\bf R}^n$:
$
\phi =\phi (x^1,\dots,x^n)(dx^1\wedge\cdots\wedge dx^n)^{\lambda }.
$
For example, in the case $\lambda=0$, 
one has the standard $\Diff(M)$-action on functions,
the value $\lambda=1$ corresponds to the $\Diff(M)$-action on
differential $n$-forms.

\vskip 0,3cm

Consider the operator of {\it Lie derivative} 
on tensor-densities of degree $\lambda$:
$$
L^{\lambda}_X(\phi)=X^i\frac{\partial\phi}{\partial x^i}+\lambda\partial_i(X^i)\phi,
$$
where $X\in\Vect({\bf R}^n)$.
(Note, that this formula does not depend on the choice of local coordinates.)
The 1-parameter family of $\Vect({\bf R}^n)$-modules 
on space ${\cal D}$
is defined by the commutator:
\begin{equation}
\ad L^{\lambda}_X(A):=L^{\lambda}_X\circ A-A\circ L^{\lambda}_X.
\label{ad}
\end{equation}

\vskip 0,3cm

{\bf Notation 2.2}. Let us denote:
\hfill\break
(a) ${\cal F}_{\lambda}$ the $\Diff({\bf R}^n)$- and $\Vect({\bf R}^n)$-modules 
of tensor-densities on ${\bf R}^n$ of degree $\lambda$.
\hfill\break
(b) ${\cal D}_{\lambda}$ and
${\cal D}^k_{\lambda}\subset{\cal D}_{\lambda}$  the $\Diff({\bf R}^n)$- and 
$\Vect({\bf R}^n)$-modules 
of differential operators defined by (\ref{act}) and (\ref{ad}).

\subsection{Main definition}

It is well known (see Sections 4.2 -- 4.4 for the precise statements), that
the notion of a symbol of a differential operator is not intrinsic.
Only the principal symbol has a geometrical meaning.
In other words,
there is no symbol map equivariant with respect to 
$\Diff({\bf R}^n)$
(or $\Vect({\bf R}^n)$).

Consider the restriction of the $\Vect({\bf R}^n)$-modules ${\cal D}_{\lambda}$ to
the projective Lie algebra $sl_{n+1}$.
Introduce the coordinates $(x^i,\xi_i),\; i=1,\dots,n$ on $T^*{\bf R}^n$,
then ${\rm Pol}(T^*{\bf R}^n)$ is the space of polynomials in $\xi$
with coefficients in $C^{\infty}({\bf R}^n)$.
We are looking for
a symbol map on ${\bf R}^n$:
$$
\sigma^{\lambda}:
{\cal D}_{\lambda}
\to
{\rm Pol}(T^*{\bf R}^n)
$$
equivariant with respect to the action of $sl_{n+1}$.

Let us formulate the first main result of this paper.

\proclaim Theorem I. (i) For every $\lambda$, there exists a unique 
$sl_{n+1}$-equivariant symbol map
$\sigma^{\lambda}$,
such that for every $A\in{\cal D}$, the higher order term
of polynomial $\sigma^{\lambda}_A(\xi)$
coincides with the principal symbol of $A$.
\hfill\break
(ii) $\sigma^{\lambda}$ maps a differential operator
$A=a_k^{i_1\dots i_k}\partial_{i_1}\cdots\partial_{i_k}$ to the polynomial
$
\sigma^{\lambda}_A(\xi)=\break
\sum_{m=0}^{k}
\bar a_{k-m}^{i_1\dots i_{k-m}}
\xi_{i_1}\cdots
\xi_{i_{k-m}}
$
with the coefficients:
$$
\bar a_{k-m}
=
\sum_{m=0}^{k}
c^k_{k-m}\;a_k^{(m)}
$$
where $a_k^{(m)}$ is the ``divergence'':
$(a_k^{(m)})^{i_1\dots i_{k-m}}=
\sum_{j}\partial_{j_1}\cdots\partial_{j_m}
(a_k^{i_1\dots i_{k-m}j_1\dots j_m})$
and $c_{k-m}^k$ are constants:
\begin{equation}
c_{k-m}^k=(-1)^m\frac{{k \choose m}{(n+1)\lambda+k-1 \choose m}}{{2k+n-m \choose m}}
\label{exp}
\end{equation}
\par

\vskip 0,3cm

{\bf Remark 2.3}.
The condition of $sl_{n+1}$-equivariance is already sufficient to determine each 
term $\bar a_{w-m}$ up to a constant.
The supplementary condition, that the higher order term
of polynomial $P^{\lambda}_A(\xi)$
coincides with the principal symbol, fixes the normalization.

\proclaim Theorem I'. The unique (up to normalization)
$sl_{n+1}$-equivariant quantization map is the map
inverse to the $sl_{n+1}$-equivariant symbol map.
It
associates to a monomial
$\bar a_k=\bar a^{i_1\dots i_k}_k\xi_{i_1}\cdots\xi_{i_k}$
the differential operator:
$$
(\sigma^{\lambda})^{-1}(\bar a_k)
=
\sum_{m=0}^{k}
\bar c_{k-m}^k
\;({\overline a_k}^{(m)})^{i_1\dots i_{k-m}}
\partial_{i_1}\cdots
\partial_{i_{k-m}}
$$
where
\begin{equation}
\bar c_{k-m}^k=\frac{{k \choose m}{(n+1)\lambda+k-1 \choose m}}{{2k+n-1 \choose m}}
\label{exp'}
\end{equation}\par

\vskip 0,3cm

{\bf Remark 2.4}. In the one-dimensional case ($n=1$)
the formul{\ae} (\ref{exp}), (\ref{exp'}) are 
valid also for the space of pseudodifferential operators.
In this case these formul{\ae} coincide with
(4.11) and (4.10) of \cite{cmz}.

\subsection{Example: second order operators and quadratic Hamiltonians}

Let us apply the general formul{\ae} (\ref{exp'}), (\ref{exp}) in the case
of a quadratic polynomials on $T^*{\bf R}^n$
and second order differential operators.

\vskip 0,3cm

(a). The $sl_{n+1}$-equivariant symbol of second order differential operator 
$A=a_2^{ij}\partial_i\partial_j+a_1^i\partial_i+a_0$ is:
$
\sigma^{\lambda} _A(\xi)
=
\bar a_2^{ij}\xi_i\xi_j+\bar a_1^i\xi_i+\bar a_0,
$
where 
$$
\matrix{  
\bar a_2^{ij}
=& 
a_2 ^{ij}\hfill \cr\noalign{\smallskip}
\bar a_1^i 
=& \displaystyle
a_1^i- 2\frac{(n+1)\lambda+1}{n+3}\partial_j(a_2^{ij})\hfill \cr\noalign{\smallskip}
\bar a_0 
=& \displaystyle
a_0 - \lambda \partial_i(a_1^i) + 
\lambda\frac{(n+1)\lambda+1}{n+2}\partial_i\partial_j(a_2^{ij}) \hfill \cr
}
$$

\vskip 0,3cm

(b). The $sl_{n+1}$-equivariant quantization map 
$(\sigma^{\lambda})^{-1}$ associates 
the second order differential operator 
$A=a_2^{ij}\partial_i\partial_j+a_1^i\partial_i+a_0$: 
$$
\matrix{  
 \displaystyle a_2^{ij}=\bar a_2^{ij}\hfill\cr\noalign{\smallskip}
 \displaystyle a_1^i=\bar a_1^i+
2\frac{(n+1)\lambda+1}{n+3}\partial_j(\bar a_2^{ij})\hfill\cr\noalign{\smallskip}
 \displaystyle a_0=\bar a_0+
\lambda\partial_i(\bar a_1^i)+
\lambda\frac{(n+1)((n+1)\lambda+1)}{(n+2)(n+3)}\partial_i\partial_j(\bar a_2^{ij})\hfill \cr
}
$$
with a polynomial
$\overline{P}(\xi)=\bar a^{ij}_2\xi_i\xi_j+\bar a_1^i\xi_i+\bar a_0$ on $T^*M$.

\vskip 0,3cm

We will consider in Section 7.3 the special case $\lambda=1/2$
and take a non-degenerate
quadratic form $H$.
Then, one has a metric $H^{-1}=g_{ij}\;dx^idx^j$ on ${\bf R}^n$.
The corresponding volume form $\sqrt g$ identifies $1/2$-densities
with functions: ${\cal F}_{1/2}\cong C^{\infty}({\bf R}^n)$.
The operator $(\sigma^{\lambda})^{-1}(H)$ in this case is 
interpreted as a Laplace-Beltrami operator.

\subsection{Structure of $sl_{n+1}$-module}

Space of polynomial on fibers functions
${\rm Pol}(T^*{\bf R}^n)
=
{\bf C}_{C^{\infty}({\bf R}^n)}[\xi_1,\dots,\xi_n]$
has a very simple $\Diff({\bf R}^n)$- and $\Vect({\bf R}^n)$-module structures.
It is isomorphic to the direct sum of modules of symmetric contravariant tensor fields:
$$
{\rm Pol}^k(T^*{\bf R}^n)
\cong
{\cal S}^0
\oplus\dots\oplus
{\cal S}^k,
$$
where ${\cal S}^l:=\Gamma (S^l(T{\bf R}^n))$
is the space of symmetric $(l,0)$-tensor fields.

\vskip 0,3cm

The structure of the $\Diff({\bf R}^n)$- and $\Vect({\bf R}^n)$-modules 
of differential operators
is much more complicated. However,
Theorem I implies that the restriction of these
modules to the projective Lie algebra $sl_{n+1}$
is, in some sense, trivial:

\proclaim Corollary 2.5.  For every $\lambda$,
${\cal D}^k_{\lambda}$
is isomorphic as a $sl_{n+1}$-module to 
${\cal S}^0
\oplus\dots\oplus
{\cal S}^k$.
In particular, $sl_{n+1}$-modules ${\cal D}_{\lambda}$
are isomorphic to each other for different values of $\lambda$.
\par

Indeed, the symbol map $\sigma^{\lambda}$ is an isomorphism
of $sl_{n+1}$-modules.

\subsection{Diagonalization of intertwining operators}

The uniqueness of $\sigma^{\lambda}$ implies
that in terms of the $sl_{n+1}$-equivariant symbol,
every isomorphisms of modules of differential operators
has a diagonal form:

\proclaim Corollary 2.6.
A linear map
$$
T:{\cal D}^k_{\lambda}\to{\cal D}^k_{\mu}
$$
is $sl_{n+1}$-equivariant if and only if the symbols
$\sigma^{\lambda} _A$ and $\sigma^{\mu} _{T(A)}$
are proportional:
$$
(\sigma^{\lambda}_0(A),\dots,\sigma^{\lambda}_k(A))
\;=\;
(\alpha_0\;\sigma^{\mu}_0(T(A)),\dots,\alpha_k\;\sigma^{\mu}_k(T(A)))
$$
where $\sigma^{\lambda}_i(A)$ 
is the homogeneous component of order $i$ of $\sigma^{\lambda}_A$,
$\alpha_0,\dots,\alpha_k\in\bf R$
are arbitrary constants.\par

In other words, the map 
$\sigma T\sigma^{-1}$
on
${\cal S}^0
\oplus\dots\oplus
{\cal S}^k$
multiplies each term of the direct sum by a constant.

\section{Polynomial language, proof of Theorem I}

In this section we prove Theorems I and I'.

Using the ``standard'' symbol, we represent 
every differential operator on ${\bf R}^n$ as polynomials
on $T^*{\bf R}^n$.
The operator (\ref{A})
corresponds to the polynomial
$$
P_A=\sum _{l=0}^ka^{i_1\dots i_l}\xi_{i_1}\cdots\xi_{i_l},
$$
where $\xi_i$ are coordinates on the fiber dual to $x^i$.

Following \cite{lec}, we rewrite the conditions of
equivariance in terms of differential equations.

\subsection{Operator of Lie derivative}

The operator $L^{\lambda}_X(f)$ 
(Lie derivative of tensor-densities of degree $\lambda$)
is given by a bilinear map
$L^{\lambda}:\Vect(M)\otimes{\cal F}_{\lambda}(M)\to{\cal F}_{\lambda}(M)$.
Thus, $L^{\lambda}$ is a differential operator
on sections of a certain vector bundle over $M\times M$.
Introducing coordinates $\xi,\eta$ on fibers of $T^*{\bf R}^n\times T^*{\bf R}^n$,
one obtains the corresponding polynomial on $T^*{\bf R}^n\times T^*{\bf R}^n$:
$$
P_{L^{\lambda}_X(f)}=
X^i\eta_if+\lambda X^i\xi_if=
\langle X,\eta\rangle f+\lambda\langle X,\xi\rangle f
$$
(the first term corresponds to the derivatives of $f$
and the second one corresponds to the derivatives of $X$).

\vskip 0,3cm

Let us give the polynomial realization
for the operator $\ad L^{\lambda}$ 
(of commutator with Lie derivative (\ref{ad})):
$\ad L^{\lambda}:\Vect(M)\otimes{\cal D}\to{\cal D}$.
The polynomial $P_{[L^{\lambda}_X,A]}$
is defined on $T^*{\bf R}^n\times T^*{\bf R}^n\times T^*{\bf R}^n$
(with coordinates $\xi,\eta,\zeta$ on fibers).
The explicit formula is as follows (see \cite{lec}):
$$
P_{[L^{\lambda}_X,A]}=
\langle X,\eta\rangle\cdot P_A
-\tau_{\xi}P_A\cdot P_X
-\lambda\tau_{\xi}P_A\cdot\langle X,\xi\rangle,
$$
where $P_A=\sum a_l^{i_1\dots i_l}\zeta_{i_1}\cdots\zeta_{i_l}$
and $P_X=X^i\zeta_i$,
coordinates $\xi$ and $\eta$
correspond to derivatives of $X$ and $A$ respectively,
$\tau_{\xi}$ acts on a polynomial $P(\zeta)$ to give
the polynomial:
$$
\tau_{\xi}P=P(\zeta+\xi)-P(\zeta).
$$
This formula is a result of easy direct computations,
the term $\tau_{\xi}(P_A)$ appears from the Leibnitz rule
by application of the differential operator $A$.

\subsection{Equivariance condition}

We are looking for a differentiable linear map
$\sigma^{\lambda}:{\cal D}\to{\rm Pol}(T^*{\bf R}^n)$
commuting with the $sl_{n+1}$-action.
Let us represent $\sigma^{\lambda}$ in a polynomial form.

Take $A\in{\cal D}$,
denote $P=P_A(\zeta)$ the corresponding polynomial.
Then, the polynomial associated with
$\sigma^{\lambda}(A)$
is defined on $T^*{\bf R}^n\times T^*{\bf R}^n$:
$$
P_{\sigma^{\lambda}(A)}=P_{\sigma^{\lambda}(A)}(\eta,P),
$$
where the coordinates $\eta$ represent derivatives of $P$.
It is linear in $P$.
Denote this polynomial ${\cal P}$.

A simple direct calculation gives
the $sl_{n+1}$-equivariance condition in terms of polynomials. 
For every $X\in sl_{n+1}$:
\begin{equation}
\matrix{
{\cal P}(\xi+\eta,\;
\langle X,\eta\rangle P-
(\lambda\langle X,\xi\rangle+P_X)\tau_{\xi}P)-
\langle X,\eta\rangle {\cal P}(\eta,P)\hfill \cr\noalign{\smallskip} 
\;\;\;\;\;\;\;\;\;\;\;\;\;\;\;\;\;\;\;\;\;\;\;\;\;\;\;\;\;\;
\;\;\;\;\;\;\;\;\;\;\;\;\;\;\;\;\;
=
(X.{\cal P})(\eta,P)-
P_X\langle\xi,\partial_\zeta\rangle{\cal P}(\eta,P)\hfill \cr
}
\label{Eq}
\end{equation}
where $(X.{\cal P})$ means the derivative of coefficients
of ${\cal P}$ in the direction of the vector field $X$.

\subsection{Proof of Theorem I} 

The Taylor expansion of
the equation (\ref{Eq}) 
with respect to $\xi$ 
is of the form:
$$
\matrix{
(X.{\cal P})(\eta,P)-
(\rho(X\otimes\xi){\cal P})(\eta,P)\hfill \cr\noalign{\smallskip} 
-{1\over2}\langle X,\eta\rangle(\xi\partial_\eta)^2{\cal P}(\eta,P)
+(\xi,\partial_\eta){\cal P}(\eta,X(\xi\partial_\zeta) P)
+{1\over2}{\cal P}(\eta,X(\xi\partial_\zeta)^2 P_Y^k)
\hfill \cr\noalign{\smallskip} 
+\lambda \langle X,\xi\rangle{\cal P}(\eta,(\xi\partial_\zeta) P)
+o(|\xi|^3)=0\hfill \cr
}
$$
where $\rho(X\otimes\xi)$ is the canonical expansion of
the natural action
$(X\otimes\xi)(Y)=\langle Y,\xi\rangle X$
of the matrix $X\otimes\xi$
on the space of polynomials such as ${\cal P}$.

\vskip 0,3cm

Vanishing of the terms of degree 0 and 1 in $\xi$
is equivalent to the invariance with respect to
the constant vector fields
and the linear vector fields respectively.

Therefore, ${\cal P}$ has constant coefficients
and it is invariant under the natural action of
$gl(n,{\bf R})$.
Hense ${\cal P}(\eta,P_Y^k)$ is a polynomial in the variables
$u:=\langle Y,\zeta\rangle$ 
and 
$v:=\langle Y,\eta\rangle$.
We denote it ${\cal P}^k(u,v)$.
It is homogeneous of degree $k$ in $u,v$.

\vskip 0,3cm

It remains to impose the invariance
under the action of the vector fields (\ref{sl})
which are of degree 2
or, equivalently, under the vector fields
$X^{\alpha}=\alpha(x)E$, where $\alpha\in{{\bf R}^n}^*$
is arbitrary.
We now consider the terms of order 2 in $\xi$. 
Since $X^{\alpha}$ is of degree 2 and since the terms of
order $\leq1$ have previously been killed,
it can be easily achieved if we note that,
for every $\theta\in{{\bf R}^n}^*$,
$\langle X^{\alpha},\theta\rangle$ is the function
$x\mapsto\alpha(x)\theta(x)$.
For instance, a term like
$\langle X,\eta\rangle(\xi\partial_{\eta})^2{\cal P}(\eta,P)$
equals
$$
2\alpha_i\eta_j\partial_{\eta_i}\partial_{\eta_j}{\cal P}(\eta,P)
$$
(where $\alpha_i,\eta_j$ are the components of
$\alpha$ and $\eta$ in the canonical basis of ${\bf R}^n$).
Taking $P=P_Y^k$ and noticing that
$$
\partial_{\eta_i}{\cal P}(\eta,P_Y^k)=
Y^i\partial_v{\cal P}^k,
$$
this may be further rewritten as
$$
2\langle Y,\alpha\rangle v\partial_v^2{\cal P}^k.
$$
In this way, one obtains the equation:
$$
u\partial_v\partial_u{\cal P}^k+
(n+k)\partial_v{\cal P}^k
+
k((n+1)\lambda +k-1){\cal P}^{k-1}=0.
$$
This leads to the relations:
$$
c_l^k
=
-\frac{k((n+1)\lambda+k-1)}{(k-l)(k+l+n)}c_l^{k-1}.
$$
where $c_l^k$ is the coefficient
of the monomial $u^lv^{k-l}$ in ${\cal P}^k$.
Together with the normalization condition $c^l_l=1$,
this is equivalent to the formula (\ref{exp}).

Hence the theorem.

\subsection{Proof of Theorem I'}

We prove Theorem I' in the same way as Theorem I.
The uniqueness of the symbol map shows that it is sufficient
to find a $sl_{n+1}$-equivariant map from 
${\rm Pol}(T^*{\bf R}^n)$ to ${\cal D}$
(satisfying the normalization condition).
Such a map is necessarily equal to $(\sigma^{\lambda})^{-1}$.
Let ${\cal P}(\eta,P)$, where $P\in{\rm Pol}(T^*{\bf R}^n)$ 
be the corresponding polynomial.

As above, one gets that the polynomial
${\cal P}(\eta,P)$ is homogeneous:
$$
{\cal P}_l(\eta,\langle Y,\zeta\rangle^k)=
\bar c_l^k\langle Y,\eta\rangle^{k-l}\langle Y,\zeta\rangle^l
$$
(the condition of equivariance with respect to the affine Lie algebra).

The equivariance with respect to the quadratic part of $sl_{n+1}$
leads to the relations:
$$
\bar c_l^k
=
\frac{(l+1)((n+1)\lambda+l)}{(k-l)(k+l+n)}\bar c_{l+1}^k.
$$
The formula (\ref{exp'}) follows.

\section{Space of differential operators on a manifold
as a module over the Lie algebra of vector fields}

Let $M$ be a smooth orientable manifold of dimension $n\geq2$ and
${\cal D}(M)$ be the space of scalar linear differential operators on $M$.
Space ${\cal D}(M)$ has a natural 1-parameter family of modules over
the Lie algebra of vector fields $\Vect(M)$ 
(see \cite{duv},\cite{lec},\cite{gar}).
The definition is analogous to those of Section 2.1:
$\Vect(M)$ acts on ${\cal D}(M)$
by the commutator with the operator of Lie derivative 
on $\lambda$-densities on $M$.

In local coordinates, the {\it invariant} formul{\ae} (\ref{A})--(\ref{ad})
are valid.

\vskip 0,3cm

In this section we will study the factor-modules
${\cal D}^k_{\lambda}(M)/{\cal D}^l_{\lambda}(M)$.
Our objective is to solve the problem of isomorphism between these modules
for different values of $\lambda$ 
and to compare them with the module
of tensor fields ${\cal S}^k\oplus\dots\oplus{\cal S}^{l+1}$.

\subsection{Action of $\Vect(M)$ on ${\cal D}^k_{\lambda}(M)$
in terms of equivariant symbols}

Let us fix an arbitrary local coordinates on $M$ and
apply the $sl_{n+1}$-equivariant symbol (\ref{exp}).

The operator:
$$
\sigma^{\lambda}\circ \ad L^{\lambda}_X\circ(\sigma^{\lambda})^{-1}:
{\cal S}^k\oplus\dots\oplus{\cal S}^0
\to
{\cal S}^k\oplus\dots\oplus{\cal S}^0
$$
expresses the action of a vector field $X$ 
on ${\cal D}^k_{\lambda}(M)$ in terms of the equivariant symbol.
This action is clearly of the form:
\begin{equation}
\matrix{
\bar a_k^X\;\;=
L_X(\bar a_k)\hfill\cr\noalign{\smallskip}
\bar a_{k-1}^X=
L_X(\bar a_{k-1})+\gamma^{\lambda}_1(X,\bar a_k)\hfill\cr\noalign{\smallskip}
\bar a_{k-2}^X=
L_X(\bar a_{k-2})+
\gamma^{\lambda}_1(X,\bar a_{k-1})+
\gamma^{\lambda}_2(X,\bar a_{k})
\hfill\cr\noalign{\smallskip}
\dots\dots\hfill\cr
}
\label{nac}
\end{equation}
where $L_X(\bar a_i)$ is the Lie derivative on the
space of symmetric contravariant tensor fields ${\cal S}^i$,
$\gamma_i^{\lambda}$ is a bilinear map:
\begin{equation}
\gamma^{\lambda}_p:\Vect(M)\otimes{\cal S}^k\to{\cal S}^{k-p}.
\label{gn}
\end{equation}
Indeed, the higher order term $\bar a_k$ is the principal symbol
and it transforms as a tensor field, the transformation law for
$\bar a_{k-1}$ has a correction term depending on $\bar a_k$ etc.

\proclaim Lemma 4.1. Operations $\gamma^{\lambda}_p$ satisfy the following properties:
\hfill\break
(a) $sl_{n+1}$-equivariance:
\begin{equation}
L_X(\gamma(Y,a))=
\gamma([X,Y],a)+
\gamma(Y,L_X(a)),\;\;\;\;\;\;\;X\in sl_{n+1},a\in{\cal S}^j;
\label{equi}
\end{equation}
(b) vanishing on $sl_{n+1}$:
\begin{equation}
\gamma(X,a)\equiv 0,\;\;\;\;\;\;\;X\in sl_{n+1}.
\label{nul}
\end{equation}
\par

{\bf Proof}: follows from the $sl_{n+1}$-equivariance of the symbol map 
$\sigma^{\lambda}$.

\subsection{Modules ${\cal D}^k_{\lambda}/{\cal D}^{k-2}_{\lambda}$}

Consider the factor-modules
${\cal D}^k_{\lambda}(M)/{\cal D}^{k-2}_{\lambda}(M)$, 
where $k\geq2$ and $\dim M\geq2$.

\proclaim Theorem 4.2. (i) All the $\Vect(M)$-modules
$
{\cal D}^k_{\lambda}(M)/{\cal D}^{k-2}_{\lambda}(M)
$
with $\lambda\not=1/2$,
are isomorphic to each other and nonisomorphic to the direct sum
of modules of tensor fields.
\hfill\break
(ii) 
The module of differential operators on half-densities is exceptional:
$$
{\cal D}^k_{1\over2}(M)/{\cal D}^{k-2}_{1\over2}(M)
\cong
{\cal S}^k\oplus{\cal S}^{k-1}.
$$\par

{\bf Proof}. 
The action  of $\Vect(M)$
on the module ${\cal D}^k_{\lambda}(M)/{\cal D}^{k-2}_{\lambda}(M)$
is given by the two first terms 
$(\bar a_k^X,\bar a_{k-1}^X)$ of the formula (\ref{nac}).
The structure of the $\Vect(M)$-module
is, therefore, determined by 
bilinear map $\gamma^{\lambda}_1(X,\bar a_k)$.

\proclaim Proposition 4.3. Bilinear map 
$\gamma^{\lambda}_1(X,\bar a_k)$ is given by:
\begin{equation}
\gamma^{\lambda}_1(X,\bar a_k)
=
(2\lambda-1)\;
\frac{k(k-1)(n+1)}{2(2k+n-1)}\;
\langle\bar\gamma(X),a_k\rangle
\label{g1}
\end{equation}
where $\bar\gamma$ is the $sl_{n+1}$-equivariant
1-cocycle (\ref{ha}) on $\Vect({\bf R})$ and
$\langle\;,\;\rangle:{\cal S}^1_2\otimes{\cal S}^k\to{\cal S}^{k-1}$
is the contraction of tensors:
$$
\langle T,a\rangle^{i_1\dots i_{k-1}}=\frac{1}{k-1}
\sum_{s=1}^{k-1}
T_{ij}^{i_s}\; a^{i_1\dots\widehat i_s\dots i_{k-1}ij}.
$$\par

{\bf Proof}. First, remark that 
$\gamma^{\lambda}_1(X,\bar a_k)$
is a second order differential operator.
Since $\gamma^{\lambda}_1$
vanishes on $sl_{n+1}$, it is of second order in $X$ (and of zero order in $a$):
$$
\gamma_1^{\lambda}(X,\bar a_k)=
\alpha \sum_{s=1}^{k-1}
\partial_i\partial_j(X^{i_s})\;
\bar a_k^{i_1\dots\widehat i_s\dots i_{k-1}ij}+
\beta\;\partial_i\partial_j(X^j)\;
\bar a_k^{i_1\dots i_{k-1}i},
$$
where $\alpha,\beta$ are some constants.

Consider the operator $A=(\sigma^{\lambda})^{-1}(\bar a_k)=
\bar a_k^{i_1\dots i_k}\partial_{i_1}\cdots\partial_{i_k}+
\cdots$.
In the commutator $[L_X^{\lambda},A]$
it is sufficient to only the terms of order zero in $\bar a_k$
and of order $\leq2$ in $X$. One has:
$$
\matrix{
[L_X^{\lambda},A] &= & 
-k\;\bar a_k^{i_1\dots i_{k-1}i}\;\partial_i(X^{i_k})\;
\partial_{i_1}\cdots\partial_{i_k}\hfill\cr\noalign{\smallskip}
&&-\bigg(\frac{k(k-1)}{2}\;\bar a_k^{i_1\dots i_{k-2}ij}\;\partial_i\partial_j(X^{i_k})+
\lambda k\;\bar a_k^{i_1\dots i_{k-1}i}\;\partial_i\partial_j(X^j)
\bigg)
\partial_{i_1}\cdots\partial_{i_{k-1}}\hfill\cr\noalign{\smallskip}
&&+\cdots\cdots\hfill\cr
}
$$

Now, to determine the constants $\alpha,\beta$,
let us apply the symbol map $\sigma^{\lambda}$
and consider only the terms of second order in $X$
and zero order in $\bar a_k$.
One easily gets:
$$
\textstyle\alpha=-\frac{k}{2}-\;c_{k-1}^k
$$
and
$$
\beta=-\lambda k-c_{k-1}^k,
$$
where $c_{k-1}^k$ is the first nontrivial constant in the
formula (\ref{exp}) of the symbol
$\sigma^{\lambda}$.
Finally,
$$
\alpha=\frac{k(n+1)}{2(2k+n-1)}(2\lambda-1)
\;\;\;\hbox{and}\;\;\;
\beta=-\frac{k(k-1)}{2k+n-1}(2\lambda-1)
$$

The formula (\ref{g1}) follows.

\vskip 0,3cm

Let us prove Theorem 4.2, Part (i).

For every $\lambda,\mu\not=1/2$,
modules
${\cal D}^k_{\lambda}(M)/{\cal D}^{k-2}_{\lambda}(M)$
and
${\cal D}^k_{\mu}(M)/{\cal D}^{k-2}_{\mu}(M)$
are isomorphic.
The isomorphism is unique (up to a constant).
The formula (\ref{g1}) implies that
this isomorphism $A^{\lambda}\mapsto A^{\mu}$,
is as follows:
\begin{equation}
\textstyle
(\bar a_k^{\mu},\;\bar a_{k-1}^{\mu})=
(\bar a_k^{\lambda},\;\frac{2\mu-1}{2\lambda-1}\bar a_{k-1}^{\lambda}).
\label{ism}
\end{equation}
An important fact is that
this isomorphism is {\it well defined}
(does not depend on the choice of coordinates).
Indeed, the locally defined map commutes with $\Vect(M)$-action
and with $\Diff(M)$-action.
Therefore, the formula 
(\ref{ism}) does not change under the coordinate transformations.

It follows from the uniqueness of the
$sl_{n+1}$-equivariant symbol that
${\cal D}^k_{\lambda}(M)/{\cal D}^{k-2}_{\lambda}(M)
\not\cong
{\cal S}^k\oplus{\cal S}^{k-1}$, if $\lambda\not=1/2$.
It is also a corollary of the
cohomological interpretation of $\gamma_1$ (cf. Section 5.4).

Theorem 4.2, Part (i) is proven.

\subsection{Exceptional case $\lambda=1/2$}

In the case $\lambda=1/2$, the term
$\gamma_1(X,\bar a_k)$ vanishes. 
The $\Vect(M)$-action (\ref{nac})
in this case is just the standard action on 
${\cal S}^k\oplus{\cal S}^{k-1}$.

Theorem 4.2 is proven.

\vskip 0,3cm

It follows from this theorem that for $\lambda\not=1/2$,
already the second term of a symbol of
a differential operator is not intrinsically defined.
However, in the exceptional case of differential operators on
$1/2$-densities, the two first terms of the
$sl_{n+1}$-equivariant symbol have geometric sense. One obtains
the following amazing remark:

\proclaim Corollary 4.4. The
$sl_{n+1}$-equivariant symbol 
defines an equivariant map:
$$
(\sigma^{1/2}_k,\sigma^{1/2}_{k-1}):
{\cal D}^k\to
{\cal S}^k\oplus{\cal S}^{k-1}.
$$\par

The geometrical reason is as follows.
The $\Vect(M)$-module ${\cal D}^k_{1/2}$ has a symmetry:
the operator of {\it conjugation}.
Every operator $A\in{\cal D}^k_{1/2}$ can be decomposed
in a sum: $A=A_0+A_1$, where 
$A_0^*=(-1)^kA_0$ and $A_1^*=(-1)^{k-1}A_1$.
Then, $\sigma^{1/2}_k$ is the principal symbol of $A_0$ and
$\sigma^{1/2}_{k-1}$ is the principal symbol of $A_1$.

\subsection{Modules ${\cal D}^k_{\lambda}/{\cal D}^l_{\lambda}$
in multi-dimensional case}

Consider the $\Vect(M)$-modules
${\cal D}^k_{\lambda}(M)/{\cal D}^l_{\lambda}(M)$
with $k-l\geq3$ in the multi-dimensional case
$(\dim M\geq2)$.
The following result shows that there is no nontrivial diffeomorphisms
between the $\Vect(M)$-modules in this case.

\proclaim Theorem 4.5. (i) Modules
${\cal D}^k_{\lambda}(M)/{\cal D}^l_{\lambda}(M)$
and
${\cal D}^k_{\mu}(M)/{\cal D}^l_{\mu}(M)$,
where $k-l\geq3$
are isomorphic if and only if $\lambda+\mu=1$.
\hfill\break
(ii) There is no isomorphism
between the modules ${\cal D}^k_{\lambda}(M)/{\cal D}^l_{\lambda}(M)$
and the module of tensor fields
${\cal S}^k\oplus\dots\oplus{\cal S}^{l+1}$. \par

{\bf Proof}.
The isomorphism in Part (i) is given by the standard conjugation
of differential operators.
This map defines a general isomorphism
$$
*:{\cal D}_{\lambda}\to{\cal D}_{1-\lambda}
$$
(cf. \cite{duv},\cite{lec}).

To prove that there is no other isomorphisms,
it is sufficient to consider the case $k-l=3$.
Indeed, module 
${\cal D}^k_{\lambda}(M)/{\cal D}^l_{\lambda}(M)$
with $k-l\geq3$
projects on
${\cal D}^k_{\lambda}(M)/{\cal D}^{k-3}_{\lambda}(M)$.

Suppose, that modules
${\cal D}^k_{\lambda}(M)/{\cal D}^{k-3}_{\lambda}(M)$
and
${\cal D}^k_{\mu}(M)/{\cal D}^{k-3}_{\mu}(M)$,
are isomorphic. It follows from the formula (\ref{ism}),
that in terms of $sl_{n+1}$-equivariant symbol, the isomorphism
is unique (up to a constant) and given in local coordinates by:
$
\textstyle
(\bar a_k^{\mu},\;\bar a_{k-1}^{\mu},\;\bar a_{k-2}^{\mu})=
(\bar a_k^{\lambda},\;\frac{2\mu-1}{2\lambda-1}\bar a_{k-1}^{\lambda}
,\;\frac{(2\mu-1)^2}{(2\lambda-1)^2}\bar a_{k-2}^{\lambda}).
$
Indeed, its restriction to the submodule 
${\cal D}^{k-1}_{\lambda}(M)/{\cal D}^{k-3}_{\lambda}(M)$
and the projection on the factor
${\cal D}^k_{\lambda}(M)/{\cal D}^{k-2}_{\lambda}(M)$
must be isomorphism.

Now, it follows from the formula (\ref{nac})
that the last formula defines an isomorphism
if and only if
$\gamma_2^{\mu}(X,\bar a_k)=
\frac{(2\mu-1)^2}{(2\lambda-1)^2}\gamma_2^{\lambda}(X,\bar a_k)$.
A direct calculation shows that if $\lambda+\mu\not=1$, then this equality is not 
satisfied (cf. the formula (\ref{gamma2}) below).

Theorem 4.5 is proven.

\vskip 0,3cm

Part (ii) of Theorem 4.5 confirms the fact (well-known ``in practice'')
that a symbol of a differential operator can not be defined in
an intrinsic way (even in the case of differential operators on $1/2$-densities).

\subsection{Modules of second order differential operators}

Consider the modules of second order differential operators
${\cal D}^2_{\lambda}$.

These modules have been classified in \cite{duv}. The result is:
{\it for every $\lambda,\not=0,{1\over2},1$,
all the $\Vect(M)$-modules
${\cal D}^2_{\lambda}(M)$
are isomorphic to each other,
modules 
${\cal D}^2_0(M)\cong{\cal D}^2_1(M)$ and ${\cal D}^2_{1\over2}(M)
$
are particular.}

For $\lambda,\mu\not=0,{1\over2},1$, there exists a unique (up to a constant) 
intertwining operator
$$
{\cal L}^2_{\lambda,\mu}:{\cal D}^2_{\lambda}(M)\to{\cal D}^2_{\mu}(M)
$$
(see \cite{duv}).

\vskip 0,3cm

(a). Let us express ${\cal L}^2_{\lambda,\mu}$ 
in terms of the $sl_{n+1}$-equivariant symbol (\ref{exp}).
It follows from Corollary 2.6, that
the map 
$(\sigma^{\mu})^{-1}\circ{\cal L}^2_{\lambda,\mu}\circ\sigma^{\lambda}$
is diagonal.

\proclaim Proposition 4.6. The
$\Vect(M)$-module isomorphism 
${\cal L}^2_{\lambda,\mu}(A)$
is defined in terms of $sl_{n+1}$-equivariant symbol by:
$$\textstyle
(\sigma_2^{\mu}(B),\;
\sigma_1^{\mu}(B),\;
\sigma_0^{\mu}(B))=
(\sigma_2^{\lambda}(A),\;
\frac{2\lambda-1}{2\mu-1}\sigma_1^{\lambda}(A),\;
\frac{\lambda(\lambda-1)}{\mu(\mu-1)}\sigma_0^{\lambda}(A)),
$$
where
$B={\cal L}^2_{\lambda,\mu}(A)$
and $\sigma_i$ are the homogeneous components of $\sigma$.\par

{\bf Proof}. This formula follows from 
the explicit expression for the
$\Vect(M)$-action (\ref{nac})
in the case of second order operators. By straightforward calculations 
(cf. Section 5.5) one has:
$$
\matrix{
\bar a_2^X=
L_X(\bar a_2)\hfill\cr\noalign{\smallskip}
\bar a_1^X=
L_X(\bar a_1)+\gamma_1(X,\bar a_1)\hfill\cr\noalign{\smallskip}
\bar a_0^X=
L_X(\bar a_0)+
\frac{\lambda(\lambda-1)}{n+2}\bar\gamma_2(X,\bar a_2)
\hfill\cr
}
$$
where 
\begin{equation}
\bar\gamma_2(X,\bar a_2)=
2\partial_i\partial_j\partial_k(X^k)\;\bar a_2^{ij}+
2\partial_j\partial_k(X^k)\;\partial_i(\bar a_2^{ij})
-
 (n+1)\;\partial_i\partial_j(X^k)\;\partial_k(\bar a_2^{ij})
\label{g2}
\end{equation}
and $\gamma_1(X,\bar a_1)$ is given by (\ref{g1}).

Proposition 4.6 is proven.

\vskip 0,3cm

(b). Let us give an intrinsic expression for the map ${\cal L}^2_{\lambda,\mu}$.

Every second order differential operator can be expressed 
(not in a unique way) as a linear combination of: 
\hfill\break
zero-order operator $\phi\mapsto f\phi$ (of multiplication by a function),
\hfill\break
first order operator $L^{\lambda}_X$ (of Lie derivative),
\hfill\break
a symmetric expression
$
[L^{\lambda}_X,L^{\lambda}_Y]_+=
L^{\lambda}_X\circ L^{\lambda}_Y+L^{\lambda}_Y\circ L^{\lambda}_X,
$
\hfill\break
where  $f\in C^{\infty}(M),X,Y,Z\in\Vect(M)$.

\proclaim Proposition 4.7. The isomorphism ${\cal L}^2_{\lambda,\mu}$
is:
$$
\begin{array}{rcl}
{\cal L}^2_{\lambda,\mu}
([L^{\lambda}_X,L^{\lambda}_Y]_+)&=&
[L^{\mu}_X,L^{\mu}_Y]_+ \hfill\\\noalign{\smallskip}
{\cal L}^2_{\lambda,\mu}
(L^{\lambda}_Z)&=&
\frac{2\lambda-1}{2\mu-1}L^{\mu}_Z\hfill\\\noalign{\smallskip}
{\cal L}^2_{\lambda,\mu}(f)&=&
\frac{\lambda(\lambda-1)}{\mu(\mu-1)}f\hfill\\
\end{array}
$$ 
\par

{\bf Proof}.
The $sl_{n+1}$-equivariant symbol of the operator of Lie derivative is:
$
\sigma^{\lambda}_{L^{\lambda}_X}=
X^i\xi_i.
$
In the same way,
$
\sigma^{\lambda}_{[L^{\lambda}_X,L^{\lambda}_Y]_+}=
\bar a_2^{ij}\xi_i\xi_j+\bar a_1^i\xi_i+\bar a_0,
$
where
\begin{equation}
\matrix{
\bar a_2^{ij}=
X^iY^j+Y^iX^j\hfill\cr\noalign{\smallskip}
\displaystyle \bar a_1^i=
\frac{2\lambda-1}{n+3}\bigg(
2(X^i\partial_j(Y^j)+Y^i\partial_j(X^j))-
(n+1)(X^j\partial_j(Y^i)+Y^j\partial_j(X^i))
\bigg)\hfill\cr\noalign{\smallskip}
\displaystyle \bar a_0=
-\frac{2\lambda(\lambda-1)}{n+2}\bigg(
X^i\partial_i\partial_j(Y^j)+X^i\partial_i\partial_j(Y^j)
\hfill\cr
\;\;\;\;\;\;\;\;\;\;\;\;\;\;\;\;\;\;\;\;\;\;\;\;\;
\;\;\;\;\;\;\;\;\;\;\;\;\;\;\;\;\;\;\;\;\;\;\;\;\;\;\;
+
\partial_i(X^i)\partial_j(Y^j)
-
(n+1)\partial_j(X^i)\partial_i(Y^j)
\bigg)\hfill\cr
}
\label{[]+}
\end{equation}
The result follows.

\vskip 0,3cm

The explicit formula for ${\cal L}^2_{\lambda,\mu}$
in terms of coefficients of differential operators 
was obtained in \cite{duv}.

\vskip 0,3cm

{\bf Remark}. The expression for ${\cal L}^2_{\lambda,\mu}$
in terms of Lie derivatives is intrinsic, but it is a
nontrivial fact 
that it does not depend on the choice of $X,Y$ and $f$
representing the same differential operator.
The expression for ${\cal L}^2_{\lambda,\mu}$
in terms of symbols is well-defined locally, but it is
a nontrivial fact that it is invariant with respect to coordinate changes.
The two facts are corollaries of the third one: the
two formul{\ae} represent the same map.

\vfill\eject

\section{Cohomology of $\Vect({\bf R}^n)$ with coefficients in
operators on tensor fields}

The relation between the $\Vect(M)$-modules of differential operators
${\cal D}_{\lambda}$
and the cohomology of $\Vect(M)$
with coefficients in
${\rm Hom}({\cal S}^k,{\cal S}^l)$
was noticed
(in the particular case of second order differential operators)
in \cite{duv}.
The measure of the difference between
$\Vect(M)$-module
${\cal D}^k_{\lambda}$ and the module of 
symmetric contravariant tensor fields
${\cal S}^0\oplus\cdots\oplus{\cal S}^k$
is represented by a class of the
first cohomology group:
$
{\rm H}^1(\Vect(M);\;{\rm Hom}({\cal S}^k(M),{\cal S}^l(M))).
$

\vskip 0,3cm

In this section we calculate
the first group of {\it differentiable} cohomology of $\Vect({\bf R}^n)$
{\it vanishing on the subalgebra} $sl_{n+1}$:
$$
{\rm H}^1(\Vect({\bf R}^n),sl_{n+1};\;{\rm Hom}_{diff}({\cal S}^k,{\cal S}^l)),
$$
where
${\rm Hom}_{diff}({\cal S}^k,{\cal S}^l))$
is the space of differential operators from
${\cal S}^k$ to ${\cal S}^l$.
This cohomology group is defined using the cochains on $\Vect({\bf R}^n)$
vanishing on $sl_{n+1}$ (cf.\cite{fuc}).
We calculate explicitly the cocycles representing
nontrivial cohomology classes and interpret them as the obstruction
for existence of equivariant symbol map.

\subsection{Modules of differential operators and cohomology}

Each term $\gamma_i^{\lambda}$ in (\ref{nac}) defines
a linear map
$c_i^{\lambda}:\Vect({\bf R}^n)\to{\rm Hom}({\cal S}^k,{\cal S}^{k-i})$
by:
$$
c_i^{\lambda}(X):=\gamma_i^{\lambda}(X,.).
$$
All the maps $c_i^{\lambda}$ are differentiable and vanish on 
the projective subalgebra $sl_{n+1}$.

The following two remarks follows from the fact
that the formula (\ref{nac}) is a 
$\Vect({\bf R}^n)$-action.

\vskip 0,3cm

(a) Operator $c_1^{\lambda}$ is a {\it 1-cocycle}:
$$
[L_X,c_1(Y)]-[L_Y,c_1(X)]-c_1([X,Y])=0,
$$
for every $X,Y\in\Vect({\bf R}^n)$.

(b)
If $\gamma_1^{\lambda}\equiv0$, that means
in the case $\lambda={1\over2}$,
it is clear that
the two maps $c_2^{1/2}$ and $c_3^{1/2}$ are 1-cocycles.
We will show that $c_3^{1/2}\equiv0$.

\vskip 0,3cm

Note, that for general values of $\lambda$, 
operator $c_2^{\lambda}$ satisfies
the relation: ${\rm d}c_2=[c_1,c_1]_+$,
where $[\;,\;]_+$ is the Massey product.

\vskip 0,3cm

{\bf Important remark}. 
The cocycles $c_1^{\lambda}$ and $c_2^{1/2}$ 
are nonzero and represent 
{\it nontrivial}
classes of cohomology.
This fact follows immediately from the uniqueness of 
$sl_{n+1}$-equivariant symbol (Theorem I).
It follows also directly from the fact that
these cocycles are $sl_{n+1}$-equivariant.
Indeed,
if $c={\rm d}b$, then $b$ is $sl_{n+1}$-equivariant,
but there is no $sl_{n+1}$-equivariant
operators from ${\cal S}^k$ to ${\cal S}^l$.

\subsection{Bilinear $sl_{n+1}$-equivariant operators}

In this section we classify the
bilinear $sl_{n+1}$-equivariant differential operators
(\ref{gn})
vanishing on the projective subalgebra
$sl_{n+1}\subset\Vect({\bf R}^n)$.

\vskip 0,3cm

This means, we consider operations $\gamma_p$
satisfying the conditions  (\ref{equi}) and (\ref{nul})
(In the contrast with the one-dimensional case,
the vanishing condition does not follow from 
$sl_{n+1}$-equivariance.)

\vskip 0,3cm

{\bf Example 5.1}: {\it The transvectants}.
Let us first recall the situation in the one-dimensional case $(M=S^1,{\bf R})$.

Bilinear $sl_2$-equivariant maps on 
tensor-densities on $\bf R$ or $S^1$
$
{\cal F}_{\lambda}\otimes{\cal F}_{\mu}\to
{\cal F}_{\lambda+\mu+m}
$
were classified by Gordan \cite{gor}.
For every $\lambda,\mu$ and $m=0,1,2,\dots$, 
there exists a remarkable operation
called the transvectant of order $m$:
\begin{equation}
J_m(\phi, \psi ) 
= 
\sum_{i+j=m} (-1)^i m!
{2\lambda+m-1 \choose i} {2\mu+m-1 \choose j}  
{\phi}^{(i)} {\psi}^{(j)} 
\label{tra}
\end{equation}
where ${\phi}^{(i)}=d^i\phi/dx^i$.
The operation (\ref{tra}) is {\it unique} (up to a constant) 
for almost all values of $\lambda$ and $\mu$.

The transvectant $J_k$
with $k\geq2$ vanishes on the projective Lie algebra $sl_2$.

\vskip 0,3cm

Let us consider the multi-dimensional case $(n\geq2)$.
The following answer is rather unexpected
and completely different from the classification 
in the one-dimensional case.

\proclaim Theorem 5.2. 
The space of
$sl_{n+1}$-equivariant
operations (\ref{gn})
vanishing on $sl_{n+1}$
is as follows:
\hfill\break
(i) 
$k>p\geq2$, 
there exist 2 independent operations.
\hfill\break
(ii)
$k=p\geq2$,
there exists a unique (up to a constant) operation.
\hfill\break
(iii)
$p=1,k\geq2$,
there exists a unique (up to a constant) operation.
\hfill\break
(iv)
There is no such operations if $k=p=1$.\par

{\bf Proof}.
Let us first consider only the 
equivariance and vanishing conditions
with respect to the affine Lie algebra
(the subalgebra of $sl_{n+1}$ generated by
the constant and linear vector fields in (\ref{sl})).
A bilinear map
$\gamma_p:\Vect({\bf R}^n)\otimes{\cal S}^k\to{\cal S}^{k-p}$,
 equivariant with respect to the affine subalgebra of $sl_{n+1}$,
is given by a homogeneous differential operator
of order $p+1$
(cf. Section 3.3).
The condition, that $\gamma_p$ vanish on the affine subalgebra
means that the expression
$\gamma_p(X,a)$ does not contain terms of order $\leq1$ in $X$.

Therefore, 
$\gamma_p$
in of the form:
\begin{equation}
\matrix{
\gamma_p(X,a)&=&
\displaystyle\sum_{u=1}^{p}\bigg(\displaystyle\sum_{s=1}^{k-p}\alpha_u\;
\partial_{j_1}\cdots\partial_{j_{u+1}}(X^{i_s})\;
\partial_{j_{u+2}}\cdots\partial_{j_{p+1}}
(a^{i_s\dots{\widehat i_s}\dots i_{k-p}j_1\dots j_{p+1}})
\hfill\cr\noalign{\smallskip}&&\;\;\;\;\;+
\beta_u\;
\partial_{j_1}\cdots\partial_{j_u}\partial_k(X^k)\;
\partial_{j_{u+1}}\cdots\partial_{j_p}(a^{i_1\dots i_{k-p}j_1\dots j_p})
\hfill\cr\noalign{\smallskip}&&\;\;\;\;\;+
\delta_u\;
\partial_{j_1}\cdots\partial_{j_u}(X^k)\;
\partial_k\partial_{j_{u+2}}\cdots\partial_{j_p}(a^{i_1\dots i_{k-p}j_1\dots j_p})
\bigg)\partial_{i_1}\otimes\cdots\otimes\partial_{i_{k-p}}
\hfill\cr
}
\label{gp}
\end{equation}
where $\alpha_u,\beta_u,\delta_u$ are constants, $u=1,\dots p$,
and $\delta_1=0$.

In invariant terms (cf. Section 1.2) this formula
reads as follows:
$$
\matrix{
\gamma_p(X,a)&=&
\displaystyle\sum_{u=1}^{p}\bigg(\alpha_u\;
\langle D^{u+1}(X),
a^{(p-u-1)}\rangle
\hfill\cr\noalign{\smallskip}&&
\;\;\;\;\;\;+
\beta_u\;
\langle D^{u}(X'),
a^{(p-u)}\rangle
\hfill\cr\noalign{\smallskip}&&
\;\;\;\;\;\;+
\delta_u\;
\langle D^u(X),
D(a^{(p-u-1)})\rangle
\bigg)
\hfill\cr
}
$$
where $D$ is the differential acting on symmetric tensor fields
(cf. Section 1.2),
$X'={\rm div}X$
and $a^{(l)}$ is the divergence of order $l$ of $a$ (cf. Section 2.2).

\vskip 0,3cm

Now, it is sufficient to impose the 
equivariance ad vanishing conditions 
with respect to the quadratic vector fields
$X=x^ix^j\partial_j$ in $sl_{n+1}$.

\vskip 0,3cm

\proclaim Lemma 5.3. The equivariance condition:
$L_X(\gamma_p(Y,a))-\gamma_p(Y,L_X(a))=\gamma_p([X,Y],a)$
is equivalent to the following recurrent system of linear equations:
\begin{equation}
\left\{
\matrix{
-u(u+2)\alpha_{u+1}+(k-p)\delta_{u+1}+(p-u)(2k+n-p+u)\alpha_u=0\hfill\cr
-u(u+1)\beta_{u+1}+(u+1)\delta_{u+1}+(p-u)(2k+n-p+u)\beta_u=0\hfill\cr
-(u^2-1)\delta_{u+1}+(p-u)(2k+n-p+u)\delta_u=0\hfill\cr
(u+1)\delta_{u+1}+(u+1)\alpha_u+(k-p+u)\delta_u+(n+1)\beta_u=0\hfill\cr
}
\right .
\label{sys}
\end{equation}
where $u=1,\dots,p$ (and $\delta_1=0$).\par

Note, that the condition of equivariance
with respect to $X$
implies: $\gamma_p(X,.)\equiv0$. Indeed, the vanishing condition reads:
$$
2\alpha_1+(n+1)\beta_1+2\delta_2=0,
$$
which coincides with the last equation of the system for $u=1$.

\vskip 0,3cm

{\bf Proof of the lemma}. It can be obtained by a 
quite complicated straightforward calculation.
Let us give here a proof based on the polynomial representation
of differential operators (cf. Section 3).

Take (without loosing generality)
$a=Y^k$, where $Y$ is a vector field.
The polynomial corresponding to
$\gamma_p(X,Y^k)$ is:
$$
\matrix{
\displaystyle P_{\gamma}&=&
\;\displaystyle P_XP_Y^{k-(p+1)}\sum_{u=1}^p
\alpha_u\langle Y,\xi\rangle^{u+1}\langle Y,\eta\rangle^{p-u}
\hfill\cr\noalign{\smallskip}&&
\displaystyle +P_Y^{k-p}\langle X,\xi\rangle\sum_{u=1}^p
\beta_u\langle Y,\xi\rangle^{u}\langle Y,\eta\rangle^{p-u}
\hfill\cr\noalign{\smallskip}&&
\displaystyle +
P_Y^{k-p}\langle X,\eta\rangle\sum_{u=1}^p
\delta_u\langle Y,\xi\rangle^{u}\langle Y,\eta\rangle^{p-u}
\hfill\cr
}
$$
where 
the variables $\xi$ and $\eta$ correspond
to derivatives of $X$ and $Y$ respectively.

Let us consider
$v_1=\langle X,\xi\rangle$,
$w_1=\langle Y,\xi\rangle$,
$v_2=\langle X,\eta\rangle$,
$w_2=\langle Y,\eta\rangle$,
$x=P_X$ and $y=P_Y$
as independent variables.
One can check that the equivariance condition leads to the following
differential equations on $P_{\gamma}$:
$$
((w_2+w_1)\partial _{w_1v_2}+y\partial _{w_1x}+y\partial _{yv_2}+(n+1)\partial _{v_1})P_{\gamma}=0
$$
and
$$
\matrix{
(-w_1\partial ^2_{w_1}+v_2\partial _{w_1v_2}+x\partial _{xw_1}+w_2\partial ^2_{w_2}+
v_1\partial _{w_1v_2}\hfill\cr
\;\;\;\;\;\;\;\;\;\;\;\;\;\;\;\;\;\;\;\;\;\;\;\;\;\;\;\;\;\;
+x\partial _{v_2y}+2w_1\partial _{w_1w_2}+2y\partial _{w_2y}+(n+1)\partial _{w_2})P_{\gamma}=0\hfill\cr
}
$$
where $\partial $ means a partial derivative.

One readily obtains the system (\ref{sys}).

Lemma 5.3 is proven.

\vskip 0,3cm

{\bf Proof of Theorem 5.2}.
(i) It is easy to see, that 
the system given by first three equations of (\ref{sys})
with the condition 
$2\alpha_1+(n+1)\beta_1+2\delta_2=0,$
is of dimension 2:
every solution is defined by $\alpha_1,\beta_1$.
The last equation is a linear combination of the other ones.
This proves Part (i) of the theorem.

\vskip 0,3cm

(ii)
In the case $k=p$, one has the following supplementary
conditions: $\alpha_i=0$, since in this case
$\gamma_p(X,a)$ is a function.
In this case the system has a unique (up to a multiple)
solution defined by the value of $\beta_1$.

Remark, that in the particular case $k=p=2$,
the operation (\ref{g2}) is a solution of the considered problem.

(iii) For $p=1,k\geq2$ every solution is proportional to
the operation (\ref{g1}).

(iv) For $k=p=1$, the system (\ref{sys}) has no solutions.

Theorem 5.2 is proven.

\subsection{Calculation of the first cohomology vanishing on $sl_{n+1}$,
multi-dimensional analogue of the Gelfand-Fuchs cocycle}

The following theorem is the second main result of this paper.

\proclaim Theorem II. In the multi-dimensional case ($n\geq2$),
$$
{\rm H}^1(\Vect({\bf R}^n),sl_{n+1};\;{\rm Hom}_{diff}({\cal S}^k,{\cal S}^l))=
\left\{
\matrix{
{\bf R},\; k-l=1,l\not=0\hfill\cr 
{\bf R},\; k-l=2\hfill\cr 
0,\;\;\hbox{otherwise}\hfill\cr 
}
\right.
$$
\par

\vskip 0,3cm

{\bf Remark 5.4}. Theorem II essentially reduces the problem of
calculation of the cohomology group
${\rm H}^1(\Vect({\bf R}^n);\;{\rm Hom}_{diff}({\cal S}^k,{\cal S}^l))$
to calculation of the cohomology group\break
${\rm H}^1(sl_{n+1};\;{\rm Hom}_{diff}({\cal S}^k,{\cal S}^l))$.

\vskip 0,3cm

{\bf Proof of the theorem}.
A 1-cocycle on $\Vect({\bf R}^n)$
vanishing on a subalgebra of $\Vect({\bf R}^n)$,
is necessarily equivariant 
with respect to this subalgebra
(cf. Section 1.2).
Let us use Theorem 5.2.
We will calculate explicitly all the
$sl_{n+1}$-equivariant cocycles on 
$\Vect({\bf R}^n)$ with values in 
${\rm Hom}_{diff}({\cal S}^k,{\cal S}^{k-1})$.
The condition of $sl_{n+1}$-equivariance guarantees, that
these cocycles are nontrivial (cf. Section 5.1).

\vskip 0,3cm

(a). If $k-l=p=1$ and $k\geq2$, there exists a unique $sl_{n+1}$-equivariant operation
$\gamma_1$
proportional to (\ref{g1}).
The corresponding linear map
$c_1$ is
a nontrivial 1-cocycle.

(b). If $k=p=1$ there is no $sl_{n+1}$ equivariant operations.

\vskip 0,3cm

(ii). 
Consider the case $p=2$. There exists a two-dimensional space
of $sl_{n+1}$ operations 
$c_2:\Vect({\bf R}^n)\to{\rm Hom}_{diff}({\cal S}^k,{\cal S}^{k-2})$
vanishing on $sl_{n+1}$.

\proclaim Proposition 5.5. There exists a unique (up to a constant)
cocycle
$c_2$ on $\Vect({\bf R}^n)$
with values in ${\rm Hom}_{diff}({\cal S}^k,{\cal S}^{k-2})$
vanishing on $sl_{n+1}$. It is defined by the bilinear map
$\gamma_2$ given by the formula (\ref{gp}) with
the coefficients:
\begin{equation}
\matrix{
\alpha_1=k-2 &
\alpha_2=\frac{1}{6}(k-2)(2k+n+1)
\hfill\cr\noalign{\smallskip}
\beta_1=1 \hfill & \beta_2=1
\hfill\cr\noalign{\smallskip}
& \delta_2=-\frac{1}{2}(2k+n-3) 
\hfill\cr
}
\label{re2}
\end{equation}
Cocycle $c_2$ is nontrivial.\par

{\bf Proof}. 
It follows from Theorem 5.2, that the space
of $sl_{n+1}$-equivariant
operations (\ref{gn}) vanishing on $sl_{n+1}$
is of dimension 2 for $p=2,k>p$.
The cocycle condition reads:
$$
L_X(\gamma_2(Y,a))+\gamma_2(Y,L_X(a))-
L_Y(\gamma_2(X,a))-\gamma_2(X,L_Y(a))=
\gamma_2([X,Y],a).
$$
This relation adds one more equation: $\beta_1=\beta_2$
to the general system (\ref{sys}).
To obtain this equation, it is sufficient to collect
the terms with 
$\partial_{i_1}\partial_{i_2}(X^j)\partial_j(\partial_m(Y^m))
a^{i_1i_2\dots i_k}$.
The unique (up to a constant)
solution of the completed system of linear equations is
given by the formula (\ref{re2}).

The corresponding map $c_2$
is indeed a
1-cocycle.
One can verify this fact by a direct calculation.

Proposition 5.5 is proven.

\vskip 0,3cm

{\bf Remark 5.6}. The cocycle $c_2$ is a {\it multi-dimensional analogue of 
the Gelfand-Fuchs cocycle}
(\ref{gf}) in the following sense. In the one-dimensional case ($n=1$), the 
corresponding cocycle on $\Vect({\bf R})$ is
given as multiplication by the
Gelfand-Fuch cocycle.

We will show in Section 5.4, that 
the operation $\gamma_2$ defined by (\ref{re2})
is proportional to
$\gamma_2^{1/2}$.

\vskip 0,3cm

{\bf 3}. Consider the case $p\geq3$.
It is easy to see that there is no solutions
of the system (\ref{sys}) satisfying the cocycle
condition. Indeed, collecting the terms with
$\partial_{i_1}\partial_{i_2}(X^j)\partial_j(\partial_m(Y^m))\allowbreak
\partial_{i_3}\cdots\partial_{i_p}(a^{i_1i_2\dots i_k})$
in the cocycle relation,
one has $\beta_1=\beta_2$.
Collecting the terms with
$\partial_{i_1}\cdots\partial_{i_p}(X^j)\partial_j(\partial_m(Y^m))
a^{i_1i_2\dots i_k}$,
one has $\beta_1=\beta_p$.
The system (\ref{sys}) together with the two new equations
has no solutions.

Theorem II is proven.

\vskip 0,3cm

\proclaim Corollary 5.7. All the bilinear
maps $\gamma_{2p+1}^{1/2}$ vanish.\par

Indeed, since $\gamma_1^{1/2}=0$,
the first nonzero map $\gamma_{2p+1}^{1/2}$, where $p\geq1$,
defines a $\Vect({\bf R}^n)$-cocycle,
which is nontrivial because of $sl_{n+1}$-equivariance.

\subsection{Calculation of $\gamma_2^{\lambda}$}

Let us now give the explicit formula for the 
bilinear map
$\gamma^{\lambda}_2(X,\bar a_{k})$
from the formula (\ref{nac}) of the $\Vect({\bf R}^n)$-action
on linear differential operators.

\proclaim Proposition 5.8. The operation 
$$
\gamma^{\lambda}_2(X,\bar a_{k})=
\frac{k(k-1)}{2(2k+n-2)}
\bar\gamma^{\lambda}_2(X,\bar a_{k}),
$$
where $\bar\gamma^{\lambda}_2(X,\bar a_{k})$
is the operation (\ref{gp}) with the coefficients:
\begin{equation}
\matrix{
\displaystyle\alpha_1=-\frac{(k-2)\bigg(2(n+1)^2\lambda(\lambda-1)+
2k^2+2kn-4k+n^2-n+2\bigg)}{2k+n-1}
\hfill\cr\noalign{\bigskip}
\alpha_2=
-(k-2)\bigg((n+1)^2\lambda(\lambda-1)+\frac{1}{3}(k^2+kn+n^2-k+n)\bigg)
\hfill\cr\noalign{\bigskip}
\displaystyle\beta_1=
\frac{(4k+n-5)(n+1)\lambda(\lambda-1)+(k-2)(k-1)}{2k+n-1}
\hfill\cr\noalign{\bigskip}
\beta_2=(4k-6)(n+1)\lambda(\lambda-1)+(k-2)n
\hfill\cr\noalign{\bigskip}
\delta_2=-(n+1)^2\lambda(\lambda-1)-(k-2)(k+n-1)
\hfill\cr
}
\label{gamma2}
\end{equation}\par

{\bf Proof}. Proof is analogous to those of Proposition 4.3.

Indeed, the operator $A=(\sigma^{\lambda})^{-1}(\bar a_k)$
is as follows:
$$
A=
\bar a_k^{i_1\dots i_k}\partial_{i_1}\cdots\partial_{i_k}+
\bar c_{k-1}^k
\partial_i(\bar a_k^{i_1\dots i_{k-1}i})\partial_{i_1}\cdots\partial_{i_{k-1}}+
\cdots,
$$
where $\bar c_{k-1}^k$
is the first coefficient in the quantization map (\ref{exp'}):
$$
\textstyle \bar c_{k-1}^k=\frac{k((n+1)\lambda+k-1)}{2k+n-1}=-c_{k-1}^k.
$$
Then,
$$
\matrix{
[L_X^{\lambda},A] &= & 
(X^i\partial_i(\bar a_k^{i_1\dots i_k})
-k\;\bar a_k^{i_1\dots i_{k-1}i}\;\partial_i(X^{i_k}))\;
\partial_{i_1}\cdots\partial_{i_k}\hfill\cr\noalign{\bigskip}
&&-
\bigg(\frac{k(k-1)}{2}\;\bar a_k^{i_1\dots i_{k-2}ij}\;\partial_i\partial_j(X^{i_k})
\hfill\cr\noalign{\smallskip}
&&+
\lambda k\;\bar a_k^{i_1\dots i_{k-1}i}\;\partial_i\partial_j(X^j)
\hfill\cr\noalign{\smallskip}
&&
+\bar c_{k-1}^k\partial_i(\bar a_k^{i_1\dots i_{k-2}ij})\partial_j(X^{i_{k-1}})
\bigg)
\partial_{i_1}\cdots\partial_{i_{k-1}}\hfill\cr\noalign{\bigskip}
&&-
\bigg(\frac{k(k-1)(k-2)}{6}\;\bar a_k^{i_1\dots i_{k-3}ijl}\;
\partial_i\partial_j\partial_l(X^{i_{k-2}})
\hfill\cr\noalign{\smallskip}
&&+
\frac{(k-1)(k-2)}{2}\lambda\; \bar a_k^{i_1\dots i_{k-2}ij}\;
\partial_i\partial_j\partial_l(X^l)\hfill\cr\noalign{\smallskip}
&&+
\frac{(k-1)(k-2)}{2}\bar c_{k-1}^k\;
\partial_i(\bar a_k^{i_1\dots i_{k-3}ijl})\;\partial_j\partial_l(X^{i_{k-2}})
\hfill\cr\noalign{\smallskip}
&&+
(k-1)\bar c_{k-1}^k\lambda\;
\partial_i(\bar a_k^{i_1\dots i_{k-2}ij})\;\partial_j\partial_l(X^l)
\bigg)
\partial_{i_1}\cdots\partial_{i_{k-2}}\hfill\cr\noalign{\smallskip}
&&+
\cdots\cdots\hfill\cr
}
$$
where $\cdots$ means the terms of order $\geq2$ in $\bar a_k$,
or $\geq4$ in $X$.

Now, one has to apply the symbol map 
$\sigma^{\lambda}$
and collect third order terms to get:
$$
\matrix{
\alpha_1=
-(\frac{k(k-2)}{2}c_{k-2}^{k-1}+
2(k-2)c_{k-2}^k+\frac{(k-1)(k-2)}{2}\bar c_{k-1}^k+
(k-2)\bar c_{k-1}^kc_{k-2}^{k-1})\hfill\cr\noalign{\smallskip}
\alpha_2=
-(\frac{k(k-1)(k-2)}{6}+
\frac{k}{2}c_{k-2}^{k-1}+
2c_{k-2}^k)
\hfill\cr\noalign{\smallskip}
\beta_1=-(
2c_{k-2}^k+k\lambda c_{k-2}^{k-1}+
(k-1)\lambda c_{k-1}^k+\bar c_{k-1}^kc_{k-2}^{k-1})\hfill\cr\noalign{\smallskip}
\beta_2=-(
\frac{k(k-2)}{2}\lambda+
kc_{k-2}^{k-1}\lambda+\frac{k}{2}c_{k-2}^{k-1}+
2c_{k-2}^k)\hfill\cr\noalign{\smallskip}
\delta_2=-(
c_{k-2}^k+\frac{k}{2}c_{k-2}^{k-1}
)\hfill\cr
}
$$
where $c_{k-2}^k$ is the second coefficient in the formula
(\ref{exp}).

Substituting the explicit expressions for
$c_{k-1}^k,\bar c_{k-1}^k,c_{k-2}^{k-1}$ and  $c_{k-2}^k$,
one obtains the formula (\ref{gamma2}).

Proposition 5.8 is proven.

\vskip 0,3cm

{\bf Remark 5.9}. Operation 
$\bar\gamma^{1/2}_2$
coincides with $\gamma_2$ from Proposition 5.4.

\vfill\eject

\section{Modules of (pseudo)differential operators on $S^1$}

Consider the space of pseudodifferential operators 
on the circle\footnote{the results of this section hold also in the complex case,
e.g. when  $M$ is the upper half-plane.}: $M=S^1$
\begin{equation}
A=\sum_{i=0}^{\infty}a_{k-i}\bigg(\frac{d}{dx}\bigg)^{k-i}
\label{pse}
\end{equation}
where $a_{k-i}\in C^{\infty}(S^1)$, $k\in{\bf R}$.

Group $\Diff(S^1)$ and Lie algebra $\Vect(S^1)$
act on the space of pseudodifferential operators
in the same was as on the space of differential operators.
Denote $\Psi{\cal D}^k_{\lambda}$ the
$\Diff(S^1)$- and $\Vect(S^1)$- modules defined by the formul{\ae}
on the space of operators (\ref{pse})
(\ref{act}) and (\ref{ad}).

We will study $\Diff(S^1)$-modules
$\Psi{\cal D}^k_{\lambda}/\Psi {\cal D}^{k-l}_{\lambda}$.
In the particular case: $k\in{\bf Z}_+, l=k+1$,
this module is just the module of differential operators on $S^1$.

\vskip 0,3cm

It follows from the uniqueness of transvectants (cf. Section 3.1) that 
in terms of $sl_2$-equivariant symbol,
the action of a vector field
on a (pseudo)differential operator 
is written via the transvectants (\ref{tra}):
\begin{equation}
\matrix{
\bar a_k ^X\hfill  
&=&
L_X^{k}(\bar a_k) \hfill \cr\noalign{\smallskip} 
\bar a_{k_1} ^X \hfill  
&=& 
L_X^{k-1}(\bar a_{k_1})  \hfill  \cr\noalign{\smallskip}
\bar a_{k-2} ^X\hfill  
&=& 
L_X^{k-2}(\bar a_{k-2}) + t_k^2(\lambda)J_3(X,\bar a_k) \hfill \cr\noalign{\smallskip}
\bar a_{k-3} ^X\hfill  
&=& 
L_X^{k-3}(\bar a_{k-3}) + t_{k-1}^2(\lambda)J_3(X,\bar a_{k-1})  
+ t_k^3(\lambda)J_4(X,\bar a_k)    \hfill  \cr
\dots\hfill  \cr
\bar a_s ^X \hfill
&=& \displaystyle L_X^{s}(\bar a_s)   
+ \sum_{i=s+2}^{k} t_{i}^{i-s}(\lambda)J_{i-s+1}(X,\bar a_i) 
\hfill  \cr
\dots\hfill  \cr
}
\label{ber}
\end{equation}
where $t_i^{i-s}(\lambda)$ are some polynomials.

\subsection{Cohomology of $\Vect({\bf R})$ vanishing on $sl_2$}

In the one-dimensional case,
symmetric tensor fields are just tensor-densities:
${\cal S}^k\cong{\cal F}_{-k}$.
The result below follows from the classification
of $sl_2$-equivariant bilinear maps on tensor-densities.

\proclaim Proposition 6.1. 
$
{\rm H}^1(\Vect({\bf R}),sl_2;\;{\rm Hom}_{diff}({\cal F}_{-k},{\cal F}_{-l}))=
\left\{
\matrix{
{\bf R},\; k-l=2\hfill\cr 
{\bf R},\; k-l=3\hfill\cr 
0,\;\;\hbox{otherwise}\hfill\cr 
}
\right.
$
\par

{\bf Proof}.
In the one-dimensional case,
all the operations (\ref{gn})
are proportional to the transvectants (\ref{tra}).
One has:
$
J_1(X,a)=L^{\lambda}_X(a)
$
and
$
J_2(X,a)=X''a.
$
The corresponding linear map
$c_2:\Vect(S^1)\to{\rm Hom}({\cal F}_{\lambda},{\cal F}_{\lambda_1}$
is a coboundary:
$c_2={\rm d}b$, where
$b\in{\rm Hom}({\cal F}_{-k},{\cal F}_{-k+1})$,
$b(a)=a'$.

The transvectants
$J_3(X,a)
=
X'''a$
(multiplication by
the Gelfand-Fuchs cocycle (\ref{gf}))
and
$$
J_4(X,a)
=
sX^{(IV)}a + 2X'''a' 
$$
correspond to nontrivial cocycles
$c_3$ and $c_4$.

For $J_p$ with $p\geq5$,
the corresponding linear maps are no more cocycles.

\subsection{Classification of $\Diff(S^1)$-modules
$\Psi{\cal D}^k_{\lambda}/\Psi {\cal D}^{k-l}_{\lambda}$}

The classification of $\Diff(S^1)$-modules
$\Psi{\cal D}^k_{\lambda}/\Psi {\cal D}^{k-l}_{\lambda}$
follows from the formula (\ref{ber}).
As in the multi-dimensional case,
zeroes of the polynomials $t_{k-s}^j(\lambda)$
corresponds to exceptional modules.

Let us formulate the result 
for general values of $k$.

\proclaim Proposition 6.2. There exists an isomorphism of $\Diff(S^1)$-modules
$$
\Psi{\cal D}^k_{\lambda}/\Psi {\cal D}^{k-l}_{\lambda}
\;\;\cong\;\;
\Psi{\cal D}^k_{\mu}/\Psi {\cal D}^{k-l}_{\mu}
$$
where $k\not=O,1/2,1,3/2,\dots$ in the following cases:
\hfill\break
(i) $l\geq2$;
\hfill\break
(ii) $l=3$, if $t^2_k(\lambda),t^2_k(\mu)\not=0$;
\hfill\break
(iii) $l=4$, if 
$\lambda,\mu$ are not the roots of polynomials
$t^2_k,t^2_{k-1},t^3_k$;
\hfill\break
(iii) $l\geq4$, if and only if $\lambda+\mu=1$.
\par

The proof is analogous to the proof of Theorems 4.2 and 4.5.
We will give the explicit formul{\ae} for the polynomials
$t^2_k,t^2_{k-1}$ and $t^3_k$ in Section 6.3.

\vskip 0,3cm

{\bf Remark: the duality}.
There exists a nondegenerate natural pairing between
spaces $\Psi {\cal D}^k/\Psi {\cal D}^l$ and 
$\Psi {\cal D}^{-l-2}/\Psi {\cal D}^{-k-2}$.
It is given by so-called {\it Adler trace} \cite{adl}: if $A\in \Psi {\cal D}^k$,
where $k\in{\bf Z}$, then
$$
\hbox{tr}(A)=\int_{S^1}a_1(x)dx.
$$

Let now $A\in\Psi {\cal D}^k/\Psi DO^l$ and 
$B\in\Psi {\cal D}^{-l-2}/\Psi {\cal D}^{-k-2}$.
Put 
$$
(A,B):=
tr(\widetilde{A}\widetilde{B}),
$$
 where
$\widetilde{A}\in\Psi {\cal D}^k,\widetilde{B}\in\Psi {\cal D}^{-l-2}$ 
are arbitrary liftings of $A$ and $B$.

Adler's trace is equivariant with respect to the action (\ref{ad}).

This means that the pairing $(\;,\;)$ is well-defined on $\Vect(S^1)$-modules.
Indeed, $([L^{\lambda}_X,A],B)+(A,[L^{\lambda}_X,B])=0$ for every
$X\in \Vect(S^1)$ (see \cite{cmz} for the details and interesting properties
of the transvectants).

\subsection{Relation to the Bernoulli polynomials}

The polynomials $t_k^j(\lambda)$ are particular cases of
the coefficients in the $SL_2$-equivariant $\star$-product
considered in \cite{cmz}.
We will not give the explicit formula here (see the formula (4.3) of \cite{cmz}).

Let us give here first examples of polynomials $t_k^j(\lambda)$.
$$
\matrix{
\displaystyle t_k^2(\lambda)=
\frac{k(k-1)}{2k-1}\bigg(\lambda^2-\lambda-\frac{(k+1)(k-2)}{12}\bigg) \hfill \cr\noalign{\smallskip}
\displaystyle t_k^3(\lambda)=
\frac{k}{6}\lambda(2\lambda-1)(\lambda-1) \hfill \cr
}
$$

Already these examples evoke an idea about the relation of  polynomials $t_k^j(\lambda)$
to the well-known {\it Bernoulli polynomials}. Indeed,
$$
t_k^2(\lambda)=
\frac{k(k-1)}{2k-1}\bigg(B_2(\lambda)-\frac{k(k-1)}{12}\bigg),\;\;\;
t_k^3(\lambda)=\frac{k}{12}B_3(\lambda)
$$
where $B_s$ is the Bernoulli polynomial of degree $s$, e.g.:
$$
\matrix{
B_0(x)=1,\;\;\;
B_1(x)=x-1/2,\;\;\;
B_2(x)=x^2-x+1/6,\hfill\cr
B_3(x)=x^3-3x^2/2+x/2,\;\;\;B_4(x)=x^4-2x^3+x^2-1/30,\hfill\cr
B_5(x)=x^5-5x^4/2+5x^3/3-x/6.\hfill\cr
}
$$
The next examples are:
$$
\matrix{
\displaystyle t_k^4(\lambda)=
\frac{k(k-1)(k-2)}{2(2k-3)(2k-5)}
\bigg(B_4(\lambda)+\frac{2k^2-6k+3}{24}B_2(\lambda) \hfill \cr\noalign{\smallskip}
\displaystyle \;\;\;\;\;\;\;\;\;\;\;\;\;\;\;\;\;\;\;\;
\;\;\;\;\;\;\;\;\;\;\;\;\;\;\;\;\;\;\;\;\;\;\;\;
-\frac{3k^4+18k^3-35k^2+8k+2}{480}\bigg) \hfill \cr\noalign{\smallskip}
\displaystyle t_k^5(\lambda)=
\frac{k(k-1)}{15(2k-7)}\bigg(B_5(\lambda)+\frac{5(k-1)(k-3)}{24}B_3(\lambda)\bigg) \hfill \cr
}
$$

\proclaim Proposition 6.3. Polynomials $t_k^{2j}(\lambda)$ are combinations
of $B_{2s}$ with $s=0,1,\dots,j$ and polynomials
$t_k^{2j+1}(\lambda)$ are combinations
of $B_{2s+1}$ with $s=0,1,\dots,j$
\par

{\bf Proof}. This statement is a simple corollary of the isomorphism 
${\cal D}_{\lambda}\cong{\cal D}_{1-\lambda}$. Indeed, with respect to the 
involution $\lambda'=1/2-\lambda$, the Bernoulli polynomials
verify the condition: $B_s(\lambda')=(-1)^sB_s(\lambda)$.

\section{Some generalizations and applications
of the $sl_{n+1}$-equivariant symbol}

Equivariance of the symbol map $\sigma^{\lambda}$
and the inverse quantization map
with respect to the action of $sl_{n+1}$
leads to their natural generalization.
The formul{\ae} (\ref{exp}) and (\ref{exp'})
are invariant under the linear-fractional
coordinate changes.
therefore, these formul{\ae} are well-defined
globally on the cotangent bundle of a locally projective manifold.

\vskip 0,3cm

Let us recall here the definition.

\subsection{Projective structures}

A {\it projective structure} on a manifold $M$ is defined by 
a projective atlas:
an atlas
with linear-fractional coordinate changes.

More precisely, a covering $(U_i)$ with a family of local diffeomorphisms
$\phi_i:U_i\to{\bf P}^n$ is called a projective atlas if the local transformations
$\phi_j\circ\phi_i^{-1}:{\bf P}^n\to{\bf P}^n$ are projective (i.e. are given
by the action of the group $PGL_{n+1}$ on ${\bf P}^n$).

Examples of locally projective
manifolds are: ${\bf R}^n,S^n,{\bf T}^n,S^l\times{\bf T}^m$ etc.

\vskip 0,3cm

A projective structure defines locally on $M$ an action of
the Lie group $SL_{n+1}$ by {\it linear-fractional transformations} and
a (locally defined) action of the Lie algebra $sl_{n+1}$
generated by vector fields (\ref{sl}),
for every system of local coordinates of a projective atlas.
This action is stable with respect to linear-fractional transformations
(the space of vector fields (\ref{sl}) is well-defined
globally on ${\bf P}^n$).

\vskip 0,3cm

{\bf Remarks 7.1}. (a) In the case of simply connected manifold $M$ 
endowed with a projective structure,
the local action of Lie algebra $sl_{n+1}$
is defined globally on $M$.

(b) All projective structures on a simply connected manifold are diffeomorphic
to each other.

(c) Any surface admits a projective structure,
the problem of existence of projective structures
for 3-dimensional manifolds is open.

\vskip 0,5cm

One has the following simple corollary of 
$sl_{n+1}$-equivariance of the symbol $\sigma^{\lambda}$.

\proclaim Corollary 7.2. Given a manifold $M$ endowed with a projective structure,
the symbol map $\sigma^{\lambda}$, given in 
arbitrary projective atlas
by the formula (\ref{exp}), is well defined globally on $M$.\par

\subsection{$SL_{n+1}$-equivariant star-products on $T^*M$}

Let us show that for every $\lambda$,
the
$sl_{n+1}$-equivariant quantization map (\ref{exp'})
defines a star-product on $T^*M$.
One obtains, therefore, a 1-parameter family
of $sl_{n+1}$-equivariant star-products.
All of them are equivalent to each other.

\vskip 0,3cm

Given a quantization map 
$\sigma^{-1}:{\rm Pol}(T^*M)\to{\cal D}(M)$,
let us introduce a new parameter $\hbar$.
For a homogeneous polynomial
$P$ of degree $k$ put:
$$
Q_{\hbar}(P)=\hbar^k\sigma^{-1}(P).
$$

Define a new associative but non-commutative operation 
of multiplication on ${\rm Pol}(T^*M)$:
\begin{equation}
F\star_{\hbar} G:=Q_{\hbar}^{-1}(Q_{\hbar}(F)\cdot Q_{\hbar}(G)).
\label{star}
\end{equation}
The corresponding algebra is isomorphic to the associative algebra
of differential operators on $M$.

The result of the operation (\ref{star}) is
a formal series in $\hbar$.
It has the following form:
$$
F\star_{\hbar}G=
FG+\sum_{k\geq1}\hbar^kC_k(F,G),
$$
where the higher order terms $C_k(F,G)$ are some differential operators.

\vskip 0,3cm

Recall, that such an operation of is called a {\it star-products}
if $C_1(F,G)$ coincides with the standard Poisson bracket
on ${\rm Pol}(T^*M)$, modulo symmetric in $F$ and $G$ terms:
$$
C_1(F,G)=\{F,G\}\;\;\;+\;\;\;
\hbox{terms symmetric in}\;(F,G)
$$
An elementary calculation shows that the associative operation
corresponding to the quantization map (\ref{exp'}) satisfies this property.

\subsection{$SL_{n+1}$-equivariant quantization of geodesic flow on $T^*S^n$}

Consider a nondegenerate 
quadratic form $H=g^{ij}\xi_i\xi_j$ 
on $T^*S^n$. We will apply the $sl_{n+1}$-equivariant
quantization map (\ref{exp'}) in the special case of $\lambda=1/2$.
This can be considered as a version
of quantization of the geodesic flow
for the corresponding metric $g=H^{-1}=g_{ij}dx^idx^j$.
Note, that 
different approaches to this problem have already been
considered (see \cite{det}, \cite{sni}, \cite{woo}).

The quantization map (\ref{exp'}) 
associates to
$H$ a  symmetric second order differential operator on
${\cal F}_{1/2}$. In coordinates of projective structure
this operator is given by the formula:
$$
A_H=
g^{ij}\frac{\partial}{\partial x^i}\frac{\partial}{\partial x^j}+
\partial_j(g^{ij})\frac{\partial}{\partial x^i}+
\frac{(n+1)}{4(n+2)}\partial_i\partial_j(g^{ij})
$$
It follows from the symmetricity, that $A_H$ is a
{\it Laplace--Beltrami} operator. This means,
$
A_H=\Delta+\Phi,
$
where $\Delta$ is the Laplace operator corresponding to the metric $g$ and
$\Phi$ is a function. 

\vskip 0,3cm

Let us
give the explicit formula for the potential $\Phi$
in the case when the coordinates of the projective structure
are {\it normal coordinates} for the metric $g$.

\vskip 0,3cm

Recall, that for every point,
there exist so-called normal coordinates in some neighborhood, such that
$$
\Gamma_{ij}^k=0,
\;\;\;{\rm and}\;\;\;
\partial_l(\Gamma_{ij}^k)={1\over3}(R^k_{li,j}+R^k_{lj,i})
$$
in this point.

\vskip 0,3cm

The condition that projective coordinates are normal
coordinates for a metric is a sort of ``compatibility condition''
for the projective structure and metric.
An example of such situation is the standard metric 
and the standard projective structure on $S^n$.

\proclaim Proposition 7.3. In normal coordinates,
\begin{equation}
A_H=\Delta-\frac{(n+1)}{12(n+2)}R
\label{geo}
\end{equation}
where $R$ is the scalar curvature.

{\bf Proof}. In normal coordinates,
$\Delta=g^{ij}\partial_i\partial_j$ (see \cite{det}).
The formula (\ref{geo}) is a corollary of
the expression: $\partial_i\partial_j(g^{ij})=-{1\over3}R$
(which can be verified by simple calculations).

\vskip 0,3cm

{\bf Remark 7.4}. Different methods of quantization
leads to the formula
$A_H=\Delta+cR$
with various values of the constant $c$.
The formula (\ref{geo}) gives a new value of $c$
different from those of \cite{det}, \cite{sni}, \cite{woo}.

\section{Discussion: two questions}

(a) The definition of the $sl_{n+1}$-equivariant symbol is purely geometric.
It would be very interesting to understand its relations with
the algebraic structure of the space of differential operators.

The first simple observation in this direction is
a remarkable (and a-priori unexpected)
coincidence between the coefficients
in the two formul{\ae} for the intertwining
operator ${\cal L}_{\lambda,\mu}^2$
(cf. Propositions 4.6 and 4.7).

The second observation is the formula (\ref{[]+})
for the $sl_{n+1}$-equivariant symbol of the
anti-commutator of two Lie derivatives:
$[L_X^{\lambda},L_Y^{\lambda}]_+=
L_X^{\lambda}\circ L_Y^{\lambda}+
L_Y^{\lambda}\circ L_X^{\lambda}$.
In the right hand side of the formula  (\ref{[]+})
one obtains two interesting operations in $X$ and $Y$.

The precise problem is as follows.
Calculate the $sl_{n+1}$-equivariant symbol 
of a symmetrized product of $k$ Lie derivatives:
$$
[L_{X_1}^{\lambda},\dots ,L_{X_k}^{\lambda}]_+:=
{\rm Sym}_{1,\dots ,k}(L_{X_1}^{\lambda}\circ\dots\circ L_{X_k}^{\lambda})
$$

(b) The second question concerns the cocycle $c_2$
on $\Vect({\bf R}^n)$ (cf Proposition 5.4).
Cocycle $c_2$ is related with the module of
differential operators on half-densities
${\cal D}^k_{1/2}$ (cf. Remark 5.7).
Note, that $c_2$ is well defined globally on
$\Vect({\bf RP}^n)$  and $\Vect(S^n)$ (cf. Section 7.1).

Cocycle $c_2$ is a multi-dimensional analogue
of the Gelfand-Fuchs cocycle in the form (\ref{gf}).
Indeed, in the one-dimensional case it is proportional to
$c_2(X)(a)=X'''a$.

The formula (\ref{re2}) is a result of calculations.
It would be interesting to study the algebraic
and geometric properties of cocycle $c_2$.
Namely, one can specify it for the Lie algebras of
Hamiltonian, contact or unimodular vector fields.

\vskip 0,5cm

{\it Acknowledgments}. It is a pleasure to acknowledge 
numerous fruitful discussions 
with C. Duval and his constant interest to this work.
We are grateful to J.-L. Brylinski,
Yu. I. Manin, P. Mattonet, E. Mourre
and O. Ogievetsky
for 
enlightening
discussions.

\vskip 1cm


\end{document}